\documentclass[journal, twocolumn, 9pt, letterpaper]{IEEEtran}

\usepackage[sort,compress]{cite}
\usepackage{graphicx}
\usepackage{subcaption}
\usepackage{booktabs}

\usepackage[vlined,ruled,norelsize,algo2e]{algorithm2e}

\usepackage[cmex10]{mathtools}
\usepackage{dsfont,amssymb,bm}

\usepackage[amsmath,thmmarks] {ntheorem}
\newtheorem {problem}       {Problem}

\newtheorem {theorem}       {Theorem}

\newtheorem{corollary}{Corollary}

\newtheorem {proposition}{Proposition}
\newtheorem {definition}{Definition}

\newtheorem {lemma}{Lemma}
\newtheorem {example}{Example}

\theoremstyle{plain}
\theorembodyfont{\normalfont}
\newtheorem {remark}{Remark}

\usepackage[T1]{fontenc}

\newcommand\BLANK{\mathfrak {E}}

\newcommand\IND{\mathds{1}}
\newcommand\PR {\mathds{P}}
\newcommand\EXP{\mathds{E}}

\newcommand\reals{\mathds{R}}
\newcommand\integers{\mathds{Z}}

\newcommand\TRANS{\intercal}

\newcommand\DEFINED{\coloneqq}

\newcommand\SUMN{\smashoperator{\sum_{n=-\infty}^\infty}}

\begin{document}

\title {Fundamental limits of remote estimation of autoregressive Markov processes under communication constraints}

\author{Jhelum Chakravorty and Aditya Mahajan%
\thanks{This work was supported in part by Fonds de recherche du Qu\'{e}bec -- Nature et technologies (FRQNT)  Team Grant PR-173396. Preliminary version of this work was presented in part in Allerton Conference on Communication, Control, and Computing, 2014, Conference on Decision and Control, 2014, IEEE Information Theory Workshop, 2015 and IEEE International Symposium on Information Theory, 2015.}
\thanks{The authors are with the Department of Electrical and Computer Engineering, McGill University, QC, Canada. Email: \texttt{jhelum.chakravorty@mail.mcgill.ca}, \texttt{aditya.mahajan@mcgill.ca}.}%
}

\maketitle

\begin{abstract}
The fundamental limits of remote estimation of autoregressive Markov processes
under communication constraints are presented. The remote estimation system
consists of a sensor and an estimator. The sensor observes a discrete-time
autoregressive Markov process driven by a symmetric and unimodal innovations
process. At each time, the sensor either transmits the current state of the
Markov process or does not transmit at all. The estimator estimates the Markov
process based on the transmitted observations. In such a system, there is a
trade-off between communication cost and estimation accuracy. Two fundamental
limits of this trade-off are characterized for infinite horizon discounted
cost and average cost setups. First, when each transmission is costly, we
characterize the minimum achievable cost of communication plus estimation
error. Second, when there is a constraint on the average number of
transmissions, we characterize the minimum achievable estimation error.
Transmission and estimation strategies that achieve these fundamental limits
are also identified. 
\end{abstract}


\section{Introduction}

\subsection{Motivation and literature overview}

In many applications such as networked control systems, sensor and
surveillance networks, and transportation networks, etc., data must be
transmitted sequentially from one node to another under a strict delay
deadline. In many of such \emph{real-time} communication systems, the
transmitter is a battery powered device that transmits over a wireless
packet-switched network; the cost of switching on the radio and transmitting a
packet is significantly more important than the size of the data packet.
Therefore, the transmitter does not transmit all the time; but when it does
transmit, the transmitted packet is as big as needed to communicate the
current source realization. In this paper, we characterize fundamental
trade-offs between the estimation error (or distortion) and the cost or
average number of transmissions in such systems.

In particular, we consider a sensor that observes a first-order autoregressive
Markov process. At each time instant, based on the current state of the
process and the history of its past decisions, the sensor determines whether
or not to transmit the current state. If the sensor does not transmit, the
receiver must estimate the state using the previously transmitted values. A
per-step distortion function measures the estimation error. We investigate two
fundamental trade-offs in this setup: (i) when there is a cost associated with
each communication, what is the minimum expected estimation error plus
communication cost; and (ii) when there is a constraint on the average number
of transmissions, what is the minimum estimation error. For both these cases,
we characterize the transmission and estimation strategies that achieve the
optimal trade-off.

Two approaches have been used in the literature to investigate real-time or
zero-delay communication. The first approach considers coding of individual
sequences~\cite{LinderLugosi:2001,WeissmanMerhav:2002,GyorgyLinderLugosi:2004,MatloubWeissman:2006};
the second approach considers coding of Markov
sources~\cite{Witsenhausen:1979,WalrandVaraiya:1983,Teneketzis:2006,MT:real-time,KaspiMerhav:2012,AsnaniWeissman:2013}.
The model presented above fits with the latter approach. In particular, it may
be viewed as real-time transmission, which is noiseless but expensive. In most
of the results in the literature, the focus has been on identifying sufficient
statistics (or information states) at the transmitter and the receiver; for
some of the models, a dynamic programming decomposition has also been derived.
However, very little is known about the solution of these dynamic programs.

The communication system described above is much simpler than the general
real-time communication setup due to the following feature: whenever the
transmitter transmits, it sends the current state to the receiver. These
transmitted events \emph{reset} the estimation error to zero. We exploit these
special features to identify an analytic solution to the dynamic program
corresponding to the above communication system.

A static (one shot) remote estimation problem was first considered
in~\cite{Marschak1954} in the context of information gathering in
organizations. The problem of optimal off line choice of measurement times was
considered in~\cite{Kushner1964}, whereas the problem of optimal online choice
of measurement times was considered in~\cite{Skafidas_Nerode1998}. The closely related
problem of event-based sampling (also called Lebesgue sampling) was considered
in~\cite{astrom_bern2002}. In addition, several variations of the remote
estimation problem have been considered in the literature. The most closely
related models
are~\cite{ImerBasar,RabiMB12_siam,XuHes2004a,LipsaMartins:2011,NayyarBasarTeneketzisVeeravalli:2013,MH2012},
which are summarized below. Other related work includes censoring sensors
\cite{Rago,App} (where a sensor takes a measurement and decides whether to
transmit it or not; in the context of sequential hypothesis testing),
estimation with measurement cost \cite{Athans,Geromel,WuAra} (where the
receiver decides when the sensor should transmit), sensor sleep scheduling
\cite{shuman_Asilomar_2006, SarkarCruz,SarkarCruz2,FedSo} (where the sensor is
allowed to sleep for a pre-specified amount of time); and event-based
communication \cite{Astrom_survey,MengChen,Shi_book} (where the sensor
transmits when a certain event takes place). We contrast our model
with~\cite{LipsaMartins:2011,NayyarBasarTeneketzisVeeravalli:2013,MH2012}
below. 

In~\cite{ImerBasar}, optimal remote estimation of i.i.d. Gaussian processes is
investigated under a constraint on the total number of transmissions. The
optimal estimation strategy is derived when the transmitter is restricted to
be of threshold-type. 

In~\cite{RabiMB12_siam}, the optimal remote estimation of a continuous-time
autoregressive Markov process driven by Brownian motion is considered under a
constraint on the number of transmissions. The optimal transmission strategy
is derived under an assumption on the structure of the optimal estimation
strategy. It is shown that the optimal transmission strategy is of a
threshold-type, where the thresholds are determined by solving a sequence of
nested optimal stopping problems.

In~\cite{XuHes2004a} optimal remote estimation of Gauss-Markov processes is
investigated when there is a cost associated with each transmission. The
optimal transmission strategy is derived when the estimation strategy is
restricted to be Kalman-like. 

In~\cite{LipsaMartins:2011,NayyarBasarTeneketzisVeeravalli:2013,MH2012},
optimal remote estimation of autoregressive Markov processes is investigated
when there is a cost associated with each transmission. It is assumed that the
autoregressive process is driven by a symmetric and unimodal noise process but
no assumption is imposed on the structure of the transmitter or the receiver.
Using different solution approaches
(~\cite{LipsaMartins:2011,NayyarBasarTeneketzisVeeravalli:2013} use
majorization theory while \cite{MH2012} uses person-by-person optimality), it
is shown that the optimal transmission strategy is threshold-based and the
optimal estimation strategy is Kalman-like (the precise form of these
strategies is stated in Theorem~\ref{thm:fin_hor}). Thus, the optimal
transmission and estimation strategies are easy to implement.

An immediate question is how to identify the optimal transmission and
estimation strategies for a given communication cost. It is shown
in~\cite{LipsaMartins:2011,NayyarBasarTeneketzisVeeravalli:2013,MH2012} that
the optimal estimation strategy does not depend on the communication cost
while the optimal transmission strategy can be computed by solving an
appropriate dynamic program. However, the dynamic programs presented
in~\cite{LipsaMartins:2011,NayyarBasarTeneketzisVeeravalli:2013,MH2012} do not
exploit the threshold structure of the optimal strategy. 

In this paper, we provide an alternative approach to identify the optimal
transmission strategies. We consider infinite horizon remote estimation
problem and show that there is no loss of optimality in restricting attention
to transmission strategies that use a time homogeneous threshold. To determine
the optimal threshold, we first provide computable expressions for the
performance of a generic threshold-based transmission strategy and then use
these expressions to identify the best threshold-based strategy. Thus, we show
that the structure of optimal strategies derived in
\cite{LipsaMartins:2011,NayyarBasarTeneketzisVeeravalli:2013,MH2012} is also
useful to compute the optimal strategy.

\subsection{Contributions}

We investigate remote estimation for two models of Markov processes---discrete state autoregressive Markov processes (Model A) and continuous state autoregressive Markov processes (Model B); both driven by symmetric and unimodal innovations process---under
two infinite horizon setups: the discounted setup with discount factor $\beta
\in (0,1)$ and the long term average setup, which we denote by $\beta = 1$ for
uniformity of notation. For both models, we consider two fundamental trade-offs:
\begin{enumerate}
  \item \emph{Costly communication}: When each transmissions costs
    $\lambda$~units, what is the minimum achievable cost of communication plus
    estimation error, which we denote by~$C^*_\beta(\lambda)$.

  \item \emph{Constrained communication}: When the average number of
    transmissions are constrained by~$\alpha \in (0,1)$, what is the minimum achievable
    estimation error, which we denote by $D^*_\beta(\alpha)$ and refer to as the
    \emph{distortion-transmission} trade-off.
\end{enumerate}

We completely characterize both trade-offs. In particular, 
\begin{itemize}
  \item In Model~A, $C^*_\beta(\lambda)$ is continuous, increasing,
    piecewise-linear, and concave in~$\lambda$ while $D^*_\beta(\alpha)$ is
    continuous, decreasing, piecewise-linear, and convex in~$\alpha$. We derive
    explicit expressions (in terms of simple matrix products) for the corner
    points of both these curves.

  \item In Model~B, $C^*_\beta(\lambda)$ is continuous, increasing, and
    concave in~$\lambda$ while $D^*_\beta(\alpha)$ is continuous, decreasing,
    and convex in~$\alpha$. We derive an algorithmic procedure to compute
    these curves by using  solutions of Fredholm integral equations of the
    second~kind. When the innovations process is Gaussian, we characterize how
    these curves scale as a function of the variance~$\sigma^2$.
\end{itemize}

We also explicitly identify transmission and estimation strategies that achieve
any point on these trade-off curves. For all cases, we show that: (i) there is
no loss of optimality in restricting attention to time-homogeneous strategies;
(ii) the optimal estimation strategy is Kalman-like; (iii) the optimal
transmission strategy is a randomized threshold-based strategy for Model~A and
is a deterministic threshold-based strategy for Model~B.

In addition,
\begin{itemize}
  \item In Model A, the optimal threshold as a function of $\lambda$
    or~$\alpha$ can be computed using a look-up table.
  \item In Model B, the optimal threshold as function of $\lambda$ or~$\alpha$
    can be computed using the solutions of Fredholm integral equations of the
    second~kind.
\end{itemize}

\subsection{Notation}

We use the following notation. $\integers$, $\integers_{\geq 0}$ and
$\integers_{>0}$ denote the set of integers, the set of non-negative integers
and the set of strictly positive integers, respectively. Similarly, $\reals$,
$\reals_{\geq 0}$ and $\reals_{>0}$ denote the set of reals, the set of
non-negative reals and the set of strictly positive reals, respectively.
Upper-case letters (e.g., $X$, $Y$) denote random variables; corresponding
lower-case letters (e.g. $x$, $y$) denote their realizations. $X_{1:t}$ is a
short hand notation for the vector $(X_1, \dots, X_t)$. Given a matrix $A$,
$A_{ij}$ denotes its $(i,j)$-th element, $A_i$ denotes its $i$-th row,
$A^\intercal$ denotes its transpose. We index the matrices by sets of the form
$\{-k, \dots, k\}$; so the indices take both positive and negative values. For
$k \in \integers_{> 0}$, $I_k$ denotes the identity matrix of dimension $k
\times k$, and $\mathbf{1}_k$ denotes $k \times 1$ vector of ones. 

$\langle v, w\rangle$ denotes the inner product between vectors $v$ and $w$,
$\PR(\cdot)$ denotes the probability of an event, $\EXP[\cdot]$ denotes the
expectation of a random variable, and $\IND\{\cdot\}$ denotes the indicator
function of a statement. We follow the convention of calling a sequence
$\{a_k\}_{k=0}^\infty$ increasing when $a_1 \le a_2 \le \cdots$. If all the
inequalities are strict, then we call the sequence strictly increasing.

\section{Model and problem formulation}

\subsection {Model}\label{sec:model}

\begin{figure}[!t]
  \centering
  \includegraphics[width=\linewidth]{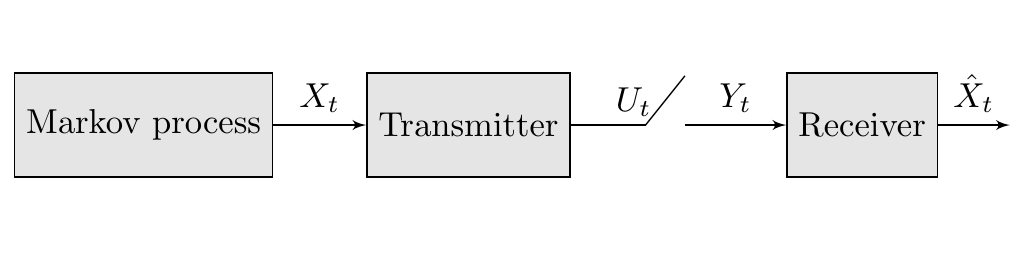}
  \caption{Block diagram of a remote estimation system.}
  \label{fig:block_diag}
\end{figure}

Consider the following two models of a discrete-time Markov process $\{X_t\}_{t=0}^\infty$ with the initial state $X_0 = 0$ and for $t \ge 0$,
\begin{equation}\label{eq:state_dyn}
  X_{t+1}=aX_t + W_t,
\end{equation}
where $\{W_t\}_{t=0}^\infty$ is an i.i.d.\ innovations process. We consider two
specific models:
\begin{itemize}
  \item \textbf{Model A:} $a, X_t, W_t \in \integers$ and $W_t$ is distributed
    according to a unimodal and symmetric pmf (probability mass function) $p$,
    i.e., for all $e \in \integers_{\ge 0}$, $p_e = p_{-e}$ and $p_e \ge
    p_{e+1}$. To avoid trivial cases, we assume $p_0$ is strictly less than 1.

  \item \textbf{Model B:} $a, X_t, W_t \in \reals$ and $W_t$ is distributed
    according to a unimodal, differentiable and symmetric pdf (probability
    density function) $\phi$, i.e., for all $e \in \reals_{\ge 0}$, $\phi(e) =
    \phi(-e)$ and for any $\delta \in \reals_{> 0}$,  $\phi(e) \ge
    \phi(e+\delta)$.
\end{itemize}

For uniformity of notation, define $\mathds X$ to be equal to  $\integers$ for
Model~A and equal to $\reals$ for Model~B. $\mathds X_{\ge 0}$ and $\mathds
X_{> 0}$ are defined similarly. 

A sensor sequentially observes the process and at each time, chooses whether
or not to transmit the current state. This decision is denoted by $U_t \in
\{0,1\}$, where $U_t = 0$ denotes no transmission and $U_1 = 1$ denotes
transmission. The decision  to transmit is made using a \emph{transmission
strategy} $f = \{f_t\}_{t=0}^\infty$, where
\begin{equation} \label{eq:transmit}
  U_t = f_t(X_{0:t}, U_{0:t-1}).
\end{equation}
We use the short-hand notation $X_{0:t}$ to denote the sequence $(X_0, \dots,
X_t)$. Similar interpretations hold for $U_{0:t-1}$. 

The transmitted symbol, which is denoted by $Y_t$, is given by 
\begin{equation*}
  Y_t = \begin{cases}
    X_t, & \text{if $U_t = 1$}; \\
    \BLANK, & \text{if $U_t = 0$},
  \end{cases}
\end{equation*}
where $Y_t = \BLANK$ denotes no transmission.

The receiver sequentially observes $\{Y_t\}_{t=0}^\infty$ and generates an estimate $\{\hat X_t\}_{t=0}^\infty$, $\hat X \in \mathds X$, using an \emph{estimation strategy} $g = \{g_t\}_{t=0}^\infty$, i.e.,
\begin{equation} \label{eq:receive}
  \hat X_t = g_t(Y_{0:t}).
\end{equation}
The fidelity of the estimation is measured by a per-step distortion
$d(X_t - \hat X_t)$.

For both models, we assume the following:
\begin{itemize}
  \item $d(0) = 0$ and for $e \neq 0$, $d(e) > 0$;
  \item $d(\cdot)$ is even, i.e., for all $e$, $d(e) = d(-e)$;
  \item $d(\cdot)$ is increasing, i.e., for $e_1 > e_2 > 0$, $d(e_1) \ge d(e_2)$;
  \item For Model~B, we assume that $d(\cdot)$ is differentiable.
\end{itemize}

We also characterize our results to the following special case of Model~B:
\begin{itemize}
  \item \textbf{Gauss-Markov model}: the density $\phi$ is zero-mean Gaussian
    with variance $\sigma^2$ and the distortion is quadratic, i.e.,
    \begin{equation*}
      \phi(e) = \frac{1}{\sqrt{2 \pi} \sigma}
        \exp{ \big(-e^2/(2 \sigma^2) \big)} \quad \text{and} \quad d(e) = e^2.
    \end{equation*}
\end{itemize}

\subsection{Performance measures}

Given a transmission and estimation strategy $(f, g)$ and a discount factor
$\beta \in (0,1]$, we define the expected distortion and the expected number
of transmissions as follows. For $\beta\in (0,1)$, the expected
\emph{discounted} distortion is given by 
\begin{equation}\label{eq:discounted_dist}
  D_\beta(f, g) \DEFINED
  (1-\beta) \EXP^{(f, g)} \Big[ \sum_{t=0}^\infty 
    \beta^t d(X_t - \hat X_t) \Bigm| X_0 = 0 \Big]
\end{equation}
and for $\beta =1$, the expected \emph{long-term average} distortion is given by
\begin{equation}\label{eq:average_dist}
  D_1(f, g) \DEFINED
  \limsup_{T \to \infty} \frac 1T \EXP^{(f, g)} \Big[ \sum_{t=0}^{T-1}
    d(X_t - \hat X_t) \Bigm| X_0 = 0 \Big].
\end{equation}

Similarly, for $\beta\in (0,1)$, the expected \emph{discounted} number of transmissions is given by 
\begin{equation}\label{eq:discounted_trans}
  N_\beta(f, g) \DEFINED
  (1-\beta) \EXP^{(f, g)} \Big[ \sum_{t=0}^\infty 
    \beta^t U_t \Bigm| X_0 = 0 \Big]
\end{equation}
and for $\beta =1$, the expected \emph{long-term average} number of transmissions is given by
\begin{equation}\label{eq:average_trans}
  N_1(f, g) \DEFINED
  \limsup_{T \to \infty} \frac 1T  \EXP^{(f, g)} \Big[ \sum_{t=0}^{T-1}
    U_t \Bigm| X_0 = 0 \Big].
\end{equation}

\begin{remark}
  We use a normalizing factor of $(1-\beta)$ to have a unified scaling for
  both discounted and long-term average setups. In particular, we will show
  that for any strategy $(f,g)$
  \begin{equation*}
    C_1(f,g;\lambda) = \lim_{\beta \uparrow 1} C_\beta (f,g;\lambda),
    \quad \text{and} \quad
    D_1(f,g) = \lim_{\beta \uparrow 1} D_\beta (f,g).
  \end{equation*}
  Similar notation is used in \cite{Altman_book}. 
\end{remark}

\subsection {Problem formulations}
We are interested in the following two optimization problems.

\begin{problem}[Costly communication]\label{prob:P1}
  In the model of Section~\ref{sec:model}, given a discount factor $\beta \in
  (0,1]$ and a communication cost $\lambda \in \reals_{>0}$, find a transmission
  and estimation strategy $(f^*,g^*)$ such that
  \begin{equation}\label{eq:opt_cost}
    C^*_\beta(\lambda) \DEFINED C_\beta(f^*,g^*;\lambda) 
    = \inf_{(f,g)} C_\beta(f,g;\lambda),
  \end{equation}
  where
  \begin{equation*}
    C_\beta(f,g;\lambda) \DEFINED D_\beta(f,g)+\lambda N_\beta(f,g)
  \end{equation*}
  is the total communication cost and the infimum in \eqref{eq:opt_cost} is
  taken over all history-dependent strategies.
\end{problem}

\begin{problem}[Constrained communication]\label{prob:P2}
  In the model of Section~\ref{sec:model}, given a discount factor $\beta \in
  (0,1]$ and a constraint $\alpha \in (0,1)$, find a transmission and estimation
  strategy $(f^*,g^*)$ such that
  \begin{equation}\label{eq:constraint}
    D^*_\beta(\alpha) \DEFINED D_\beta(f^*,g^*) 
    = \inf_{(f,g): N_\beta(f,g) \le \alpha} D_\beta(f,g),
  \end{equation}
  where the infimum is taken over all history-dependent strategies.
\end{problem}

\begin{remark}\label{rem:D_star_inf}
  It can be shown for $|a|\ge 1$ that $\lim_{\alpha \rightarrow 0}
  D^*_1(\alpha) = \infty$\footnote{For $|a|\ge 1$, a symmetric Markov chain
   as given by \eqref{eq:state_dyn} does not have a stationary distribution.
   Therefore, in the limit of no transmission, the expected long-term average
 distortion diverges to $\infty$.} and $\lim_{\alpha \rightarrow 1}
 D^*_\beta(\alpha) = 0$. 
\end{remark}

The function $D_\beta^*(\alpha)$, $\beta \in (0, 1]$, represents the minimum
expected distortion that can be achieved when the expected number of
transmissions are less than or equal to $\alpha$. It is analogous to the
distortion-rate function in Information Theory; for that reason, we call it
the \emph{distortion-transmission function}.

\section{The main results}\label{sec:main_results}

\subsection{Structure of optimal strategies}

To completely characterize the functions $C^*_\beta(\lambda)$ and
$D^*_\beta(\alpha)$, we first establish the structure of optimal transmitter
and receiver. 

\begin{theorem}[Structural results]\label{thm:inf_hor}
  Consider Problem~\ref{prob:P1} for $\beta \in (0,1]$. Then, for both
  Models~A and~B, we have the following. 
  \begin{enumerate}
    \item \emph{Structure of optimal estimation strategy}: The optimal
      estimation strategy $\hat{X}_0 = 0$ and for $t>0$ is as follows:
      \begin{align*}
        \hat{X}_t &= \begin{cases}
                     Y_t,& \mbox{if $Y_t \neq \BLANK$}\\
                     a \hat{X}_{t-1}, & \mbox{if $Y_t = \BLANK$},
                      \end{cases}
        \shortintertext{or equivalently,}
        \hat{X}_t &= \begin{cases}
                     X_t, & \mbox{if $U_t = 1$}\\
                     a \hat{X}_{t-1}, & \mbox{if $U_t \neq 1$}.
                    \end{cases}
      \end{align*}
      We denote this strategy by $g^*$.
    \item \emph{Structure of optimal transmission strategy}: Define $E_t
      \DEFINED X_t - a \hat{X}_{t-1}$, which we call the \emph{error process}.
      Then there exists a time-invariant threshold $k$ such that the
      transmission strategy 
      \begin{align}
        U_t = f^{(k)}(E_t) &\DEFINED \begin{cases}\label{eq:opt_trans}
                     1, & \mbox{if $|E_t| \ge k$}\\
                     0, & \mbox{if $|E_t| < k$}
                      \end{cases}
      \end{align}
      is optimal.
  \end{enumerate}
\end{theorem}

The proof of the theorem is given in Section~\ref{sec:proof_struct}.

Similar structural results were established for the finite horizon setup in
\cite{LipsaMartins:2011,NayyarBasarTeneketzisVeeravalli:2013,MH2012}, which we
use to establish Theorem~\ref{thm:inf_hor}. See Section~\ref{sec:proof_struct}
for details. The transmission strategy of the form~\eqref{eq:opt_trans} are
also called \emph{event-driven transmission} or \emph{delta sampling}.

\begin{remark}\label{rem:reset}
  Each transmission \emph{resets} the state of the error process to $w \in
  \mathds X$ with probability $p_w$ in Model~A and with probability density
  $\phi(w)$ in Model~B. In between the transmission, the error process evolves
  in a Markovian manner. Thus $\{E_t\}_{t=0}^\infty$ is a regenerative
  process.
\end{remark}

\subsection{Performance of generic threshold-based strategies}
\label{sec:performance_thresh}

Let $\mathcal{F}$ denote the class of all time-homogeneous threshold-based
strategies of the form~\eqref{eq:opt_trans}. For $\beta \in (0,1]$ and $e \in
\mathds X$, define the following for a system that starts in state $e$ and
follows strategy $f^{(k)}$: 
\begin{itemize}
  \item $L^{(k)}_\beta (e)$: the expected distortion until the first transmission;
  \item $M^{(k)}_\beta (e)$: the expected time until the first transmission;
  \item $D^{(k)}_\beta (e)$: the expected distortion;
  \item $N^{(k)}_\beta (e)$: the expected number of transmissions;
  \item $C^{(k)}_\beta (e;\lambda)$: the expected total cost, i.e.,
    \begin{equation*}
      C^{(k)}_\beta (e;\lambda) 
      = D^{(k)}_\beta (e) + \lambda N^{(k)}_\beta (e), \quad \lambda \ge 0.
    \end{equation*}
\end{itemize}

Note that $D^{(k)}_\beta (0) = D_\beta(f^{(k)},g^*)$, $N^{(k)}_\beta (0) =
N_\beta(f^{(k)},g^*)$ and $C^{(k)}_\beta (0;\lambda) =
C_\beta(f^{(k)},g^*;\lambda)$.

Define $S^{(k)}$ as follows:
\begin{equation*}
  S^{(k)} \DEFINED 
  \begin{cases}
    \{-(k-1), \cdots, k-1\}, & \mbox{for Model A};\\
    (-k,k),& \mbox{for Model B}.
  \end{cases}
\end{equation*}

Under strategy $f^{(k)}$, the transmitter does not transmit if $E_t \in
S^{(k)}$. For that reason, we call $S^{(k)}$ the \emph{silent set}. Define
linear operator $\mathcal{B}^{(k)}$ as follows:
\begin{itemize}
  \item \textbf{Model A}: For any $v^{(k)}: S^{(k)} \to \reals$, define
    operator $\mathcal{B}^{(k)}$ as
    \begin{equation*}
      [\mathcal{B}^{(k)}v](e) 
      \DEFINED \sum_{n \in S^{(k)}} p_{n-ae} v(n), \quad \forall e \in S^{(k)}.
    \end{equation*}
  \item \textbf{Model B}: For any $v^{(k)}: S^{(k)} \to \reals$, define
    operator $\mathcal{B}^{(k)}$ as
    \begin{equation*}
      [\mathcal{B}^{(k)}v](e)
      \DEFINED \int_{S^{(k)}} \phi(n-ae) v(n)dn, \quad \forall e \in S^{(k)}.
    \end{equation*}
\end{itemize}

Recall from Remark~\ref{rem:reset} that the state $E_t$ evolves in a Markovian
manner until the first transmission. We may equivalently consider the Markov
process until it is absorbed in $(-\infty,-k]\cup [k,\infty)$. Thus, from
balance equation for Markov processes, we have for all $e \in S^{(k)}$, 
\begin{align}\label{eq:L}
  L^{(k)}_\beta (e) &=d(e) + \beta [\mathcal{B}^{(k)} L^{(k)}_\beta](e),\\ \label{eq:M}
  M^{(k)}_\beta (e) &=1 + \beta [\mathcal{B}^{(k)} M^{(k)}_\beta](e).
\end{align}

\begin{lemma}\label{lem:propLM}
  For any $\beta \in (0,1]$, equations~\eqref{eq:L} and~\eqref{eq:M} have
  unique and bounded solutions $L^{(k)}_\beta$ and $M^{(k)}_\beta$ that are
  \begin{enumerate}
    \item[(a)] strictly increasing in $k$,
    \item[(b)] continuous and differentiable in $k$ for Model B,
    \item[(c)] $\displaystyle \lim_{\beta \uparrow 1} L^{(k)}_\beta(e) =
      L^{(k)}_1(e)$, $\displaystyle \lim_{\beta \uparrow 1} M^{(k)}_\beta(e) =
      M^{(k)}_1(e)$, for all $e$.
  \end{enumerate}
\end{lemma}
The proof of the lemma is given in Appendix~\ref{lem:proof_propLM}.

\begin{theorem}[Renewal relationships]\label{thm:renewal_DNC}
  For any $\beta \in (0,1]$, the performance of strategy $f^{(k)}$ in both
  Models~A and~B is given as follows:
  \begin{enumerate}
    \item  $D_\beta(f^{(0)}, g^*) = 0$, 
        $N_\beta(f^{(0)}, g^*) = 1$,
      and 
        $C_\beta(f^{(0)}, g^*;\lambda) = \lambda$.
    \item For $k \in \mathds X_{> 0}$,
      \begin{align*}
        D_\beta(f^{(k)}, g^*) &= \frac {L^{(k)}_\beta(0)}{M^{(k)}_\beta(0)},\\
        N_\beta(f^{(k)}, g^*) &= \frac {1}{M^{(k)}_\beta(0)} - (1-\beta),\\
      \shortintertext{and}
       C_\beta(f^{(k)}, g^*;\lambda) &= 
        \frac{ L^{(k)}_\beta(0) + \lambda }{ M^{(k)}_\beta(0)} - \lambda (1-\beta).
      \end{align*}
  \end{enumerate}
\end{theorem}
The proof of the Theorem is given in Section~\ref{sec:proof_thm_renewal}.

\begin{remark}
  There is a $- 1/(1-\beta)$ term in the expression of $N^{(k)}_\beta(0)$
  because for $k > 0$, $U_0 = 0$. Had we defined $U_0 = 1$, then we would have
  obtained the usual renewal relationship of $N^{(k)}_\beta(0) =
  1/M^{(k)}_\beta(0)$.
\end{remark}

Thus, to compute $D_\beta(f^{(k)}, g^*)$ and $N_\beta(f^{(k)}, g^*)$, one
needs to compute only
$L^{(k)}_\beta(0)$ and $M^{(k)}_\beta(0)$. Computation of the latter
expressions is given in the next section.

\begin{proposition}\label{prop:prelim_c_star_k_star}
  For both Models~A and~B,
  \begin{enumerate}
    \item $C^{(k)}_\beta(0;\lambda)$ is submodular in $(k,\lambda)$, i.e., for
      $l > k$, $C^{(l)}_\beta(0;\lambda) - C^{(k)}_\beta(0;\lambda)$ is
      decreasing in $\lambda$.

    \item Let $k^*_\beta(\lambda) = \arg\inf_{k \ge 0}
      C^{(k)}_\beta(0;\lambda)$ be the optimal $k$ for a fixed $\lambda$. Then
      $k^*_\beta(\lambda)$ is increasing in $\lambda$.
  \end{enumerate}
\end{proposition}
The proof of the proposition is in Appendix~\ref{app:prelim_c_star_k_star}.

\subsection{Computation of $L^{(k)}_\beta$ and $M^{(k)}_\beta$}\label{sec:Comp_LM}

\subsubsection{Model~A}

For Model A, the values of $L^{(k)}_\beta$ and $M^{(k)}_\beta$ can be computed
by observing that the operator $\mathcal{B}^{(k)}$ is equivalent to a matrix
multiplication. In particular, define the matrix $P^{(k)}$ as
\begin{equation*}
  P^{(k)}_{ij} \DEFINED p_{i-j}, \quad \forall i,j \in S^{(k)}.
\end{equation*}
Then, 
\begin{equation}\label{eq:withP}
  [\mathcal{B}^{(k)}v](e) = \sum_{n \in S^{(k)}} p_{n-ae} v(n)= \sum_{n \in S^{(k)}} P^{(k)}_{n,ae} v(n) = [P^{(k)} v]_{ae}.
\end{equation}

With a slight abuse of notation, we are using $v$ both as a function and a
vector. Define the matrix $Q^{(k)}$ and the vector $d^{(k)}$ as follows:
\begin{equation*}
  Q^{(k)}_\beta \DEFINED [I_{2k - 1} - \beta P^{(k)} ]^{-1}, \quad
  d^{(k)} \DEFINED [ d(-k+1), \dots, d(k-1) ]^\TRANS.
\end{equation*}
Then, \eqref{eq:L}, \eqref{eq:M} and~\eqref{eq:withP} imply the following:
\begin{proposition}\label{prop:expressLM_A}
  In Model A, for any $\beta \in (0,1]$,
  \begin{align}\label{eq:L_matrix_form}
    L^{(k)}_\beta &= [I_{2k-1} - \beta P^{(k)}]^{-1} d^{(k)}\\ \label{eq:M_matrix_form}
    M^{(k)}_\beta &= [I_{2k-1} - \beta P^{(k)}]^{-1} \mathbf{1}_{2k-1}.
  \end{align}
\end{proposition}
See Section~\ref{sec:example_A} for an example of these calculations.

\subsubsection{Model~B}
For Model B, for any $\beta \in (0,1]$, \eqref{eq:L} and~\eqref{eq:M} are
Fredholm integral equations of second kind \cite{Polyanin}. The solution can
be computed by identifying the inverse operator
\begin{equation*}
  \mathcal{Q}^{(k)}_\beta = [I-\beta \mathcal{B}^{(k)}]^{-1},
\end{equation*}
which is given by
\begin{equation}\label{eq:resolvent}
  [\mathcal{Q}^{(k)}_\beta v](e) = \int_{-k}^k R^{(k)}_\beta(e,w;a) v(w) dw,
\end{equation}
where for any given $a$, $R^{(k)}_\beta(\cdot,\cdot;a)$ is the resolvent of
$\phi$ and can be computed using the  Liouville-Neumann series.
See~\cite{Polyanin} for details. Since $\phi$ is smooth,  \eqref{eq:L}
and~\eqref{eq:M} can also be solved by discretizing the integral equation
using quadrature methods. A Matlab implementation of this approach is
available in~\cite{Atkinson08solvingfredholm}.

\subsection{Main results for Model~A}\label{sec:results_Model_A}

\subsubsection{Results for costly communication}\label{subsec:costly_A}

\begin{figure}[t]
  \subcaptionbox{}{\includegraphics[scale=0.9]{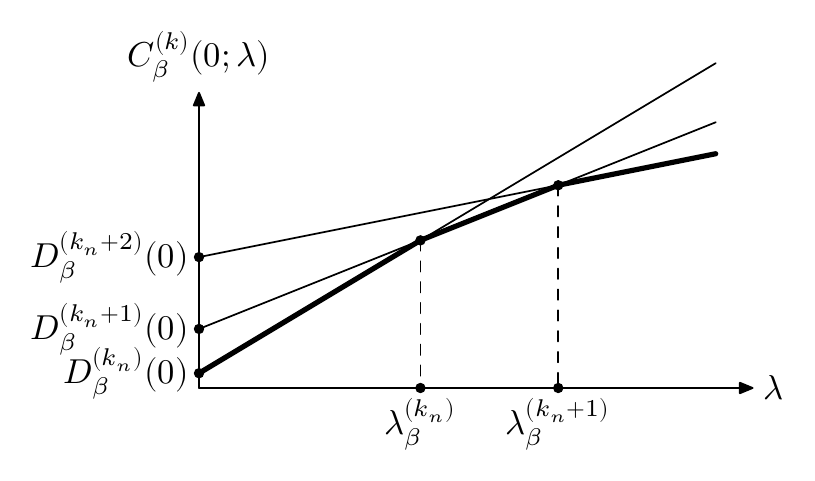}}
  \\
  % Align the y-axis of the two figures
  \null\hskip 3.5em
  \subcaptionbox{}{\includegraphics[scale=0.9]{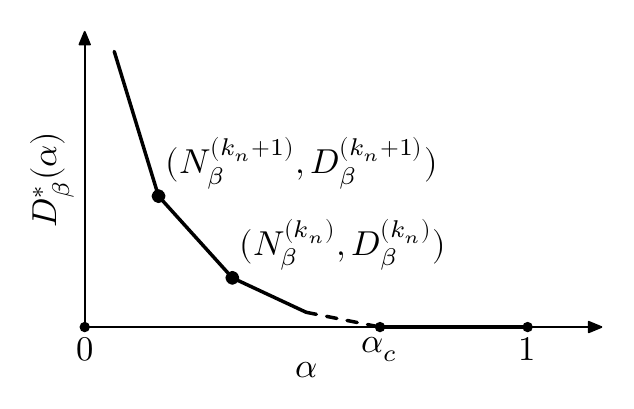}}
  \caption{In Model~A, (a)~the optimal costly communication cost
    $C^*_\beta(\lambda)$; (b)~the distortion-transmission function
    $D^*_\beta(\alpha)$.} 
  \label{fig:modelA}
\end{figure}

\begin{theorem}\label{thm:DIS_AVG}
  For $\beta \in (0, 1]$, let $\mathds K$ denote $\{ k \in \integers_{\ge
  0} : D^{(k+1)}_\beta(0) > D^{(k)}_\beta(0)\}$. For $k_n \in \mathds K$, define:
  \begin{equation}\label{eq:lambda}
    \lambda^{(k_n)}_\beta \DEFINED 
    \frac {D^{(k_{n+1})}_\beta(0) - D^{(k_n)}_\beta(0) }
          { N^{(k_n)}_\beta(0) - N^{(k_{n+1})}_\beta(0) }.
  \end{equation}
  Then, we have the following.
  \begin{enumerate}
    \item For any $k_n \in \mathds K$ and any $\lambda \in
      (\lambda^{(k_{n-1})}_\beta, \lambda^{(k_n)}_\beta]$, the strategy $f^{(k_n)}$ is
      optimal for Problem~\ref{prob:P1} with communication cost $\lambda$.

    \item The optimal performance $C^*_\beta(\lambda)$ is continuous, concave,
      increasing and piecewise linear in $\lambda$. The corner points of
      $C^*_\beta(\lambda)$ are given by $\{(\lambda^{(k_n)}_\beta,
      D^{(k_n)}_\beta(0)+\lambda^{(k_n)}_\beta N^{(k_n)}_\beta(0))\}_{k_n \in
      \mathds K}$
      (see Fig~\ref{fig:modelA}(a)).
  \end{enumerate}
\end{theorem}
The proof of the theorem is given in Section~\ref{sec:proofs_A}.

\subsubsection{Results for constrained communication}\label{subsec:constrained_A}

To describe the solution of Problem~\ref{prob:P2}, we first define Bernoulli randomized strategy and Bernoulli randomized simple strategy~\cite{Sennott_dis}.

\begin{definition}
  Suppose we are given two (non-randomized) time-homogeneous strategies $f_1$
  and $f_2$ and a randomization parameter $\theta \in (0,1)$. The
  \emph{Bernoulli randomized strategy} $(f_1, f_2, \theta)$ is a strategy that
  randomizes between $f_1$ and $f_2$ at each stage; choosing $f_1$ with
  probability $\theta$ and $f_2$ with probability $(1-\theta)$. Such a
  strategy is called a Bernoulli randomized \emph{simple} strategy if $f_1$
  and $f_2$ differ on exactly one state, i.e.,\ there exists a state $e_0$
  such that
  \begin{equation*}
    f_1(e) = f_2(e), \quad \forall e \neq e_0.
  \end{equation*}
\end{definition}

\begin{theorem}\label{thm:constrained_f_star_D}
  For any $\beta \in (0,1]$ and $\alpha \in (0,1)$, define
  \begin{align}
    k^*_\beta(\alpha) 
    &= \sup \{k \in \integers_{\geq 0}: N_{\beta}(f^{(k)}, g^*) \geq \alpha\}
    \notag\\
    &= \sup \Big\{k \in \integers_{\geq 0}: M^{(k)}_\beta \leq \frac{1}{1+\alpha -\beta}\Big\}
    \label{eq:k_star} \\
    \shortintertext{and}
    \theta^*_\beta(\alpha) 
    &= \frac{\alpha -  N_{\beta}(f^{(k^*_\beta(\alpha)+1)}, g^*)}
        { N_{\beta}(f^{(k^*_\beta(\alpha))}, g^*) 
        - N_{\beta}(f^{(k^*_\beta(\alpha)+1)}, g^*)} 
        \notag \\
    &= \frac{M^{(k^*+1)}_\beta - \frac{1}{1+\alpha-\beta}}
            {M^{(k^*+1)}_\beta - M^{(k^*)}_\beta}.
        \label{eq:theta_star}
  \end{align}
  For ease of notation, we use $k^*=k^*_\beta(\alpha)$ and $\theta^* =
  \theta^*_\beta(\alpha)$. 

  Let $f^*$ be the Bernoulli randomized simple strategy $(f^{(k^*)},
  f^{(k^*+1)}, \theta^*)$, i.e.,
  \begin{equation}\label{eq:f_star}
    f^*(e) = \begin{cases}
      0,& \text{if  $|e|<k^*$}; \\
      0, & \text{w.p. $1-\theta^*$, if $|e|=k^*$};\\
      1, & \text{w.p. $\theta^*$, if $|e|=k^*$};\\
      1,& \text{if $|e|>k^*$}.
    \end{cases}
  \end{equation}
  Then 
  \begin{enumerate}
    \item $(f^*, g^*)$ is optimal for the constrained Problem~\ref{prob:P2}
      with constraint $\alpha$. 

    \item Let $\alpha^{(k)} = N_\beta(f^{(k)}, g^*)$. Then, for $\alpha \in
      (\alpha^{(k+1)}, \alpha^{(k)})$, $k^*=k$ and $\theta^* = (\alpha -
      \alpha^{(k+1)})/ (\alpha^{(k)} - \alpha^{(k+1)})$, and the
      distortion-transmission function is given by
      \begin{equation}\label{eq:D_opt_rand}
        D^*_\beta(\alpha) = \theta^* D_\beta^{(k)} + (1-\theta^*) D_\beta^{(k+1)}.
      \end{equation}
      Moreover, the distortion-transmission function is is continuous, convex,
      decreasing and piecewise linear in $\alpha$. Thus, the corner points of
      $D^*_\beta(\alpha)$ are given by $\{(N^{(k)}_\beta(0),
      D^{(k)}_\beta(0))\}_{k=1}^\infty$ (see
      Fig~\ref{fig:modelA}(b)).
  \end{enumerate}
\end{theorem}
The proof of the theorem is given in Section~\ref{sec:proofs_A}.

\begin{corollary}
  In Model~A, for any $\beta \in (0,1]$,
  \begin{equation*}
    D_\beta(f^{(1)},g^*) = 0,
    \quad \text{and} \quad 
    N_\beta(f^{(1)},g^*) = \beta (1-p_0)\coloneqq \alpha_c.
  \end{equation*}
\end{corollary}

\subsection{Main results for Model~B}\label{sec:results_Model_B}

\subsubsection{Results for costly communication}\label{subsec:costly_B}
Let $\partial_k D^{(k)}_\beta $, $\partial_k N^{(k)}_\beta$ and $\partial_k
C^{(k)}_\beta$ denote the derivative of $D^{(k)}_\beta$, $N^{(k)}_\beta$ and
$C^{(k)}_\beta$ with respect to $k$ (in Lemma~\ref{lem:diff_DNC} we show
that $D^{(k)}_\beta$, $N^{(k)}_\beta$ and $C^{(k)}_\beta$ are differentiable
in~$k$).

\begin{theorem}\label{thm:lambda_k}
  For $\beta \in (0,1]$, we have the following.
  \begin{enumerate}
    \item If the pair $(\lambda, k)$ satisfies the following
      \begin{equation}\label{eq:lambda_k}
        \lambda = -\frac{\partial_k D^{(k)}_\beta(0)}{\partial_k N^{(k)}_\beta(0)},
      \end{equation} 
      then, the strategy $(f^{(k)}, g^*)$ is optimal for Problem~\ref{prob:P1}
      with communication cost $\lambda$. Furthermore, for any $k>0$, there
      exists a $\lambda \ge 0$ that satisfies~\eqref{eq:lambda_k}.
    \item The optimal performance $C^*_\beta(\lambda)$ is continuous, concave
      and increasing function of $\lambda$. 
  \end{enumerate}
\end{theorem}
The proof of the theorem is given in Section~\ref{sec:proofs_B}.
Algorithm~\ref{algo:C_star} shows how to compute $C^*_\beta(\lambda)$.

\begin{algorithm2e}[!tb]
  \SetKwInOut{Input}{input}
  \SetKwInOut{Output}{output}
  \SetKwFor{ForAll}{forall}{do}{}
  \DontPrintSemicolon
  \Input{$\lambda \in \reals_{>0}$, $\beta \in (0,1]$, $\varepsilon \in \reals_{>0}$}
  \Output{$C^{(k^\circ)}_\beta(\lambda)$, where $|k^\circ - k^*_\beta(\lambda)| < \varepsilon$}
  Let $\lambda^*_\beta(k)$ denote the left-hand side of~\eqref{eq:lambda_k}\;
  Pick $\underline k$ and $\bar k$ such that $\lambda^*_\beta(\underline k) < \lambda < \lambda^*_\beta(\bar k)$\;
  $k^\circ = (\underline k + \bar k)/2$\;
  \While{$|\lambda^*_\beta(k^\circ) - \lambda| > \varepsilon$}{
    \uIf{$\lambda^*(k^\circ) < \lambda$}{$\underline k= k^\circ$}
    \Else{$\bar k = k^\circ$}
    $k^\circ = (\underline k + \bar k)/2$\;
  }
  \Return{$D^{(k^\circ)}_\beta(0) + \lambda N^{(k^\circ)}_\beta(0)$\; }
  \caption{Computation of $C^*_\beta(\lambda)$}\label{algo:C_star}
\end{algorithm2e}

\subsubsection{Results for constrained communication}\label{subsec:constrained_B}
\begin{theorem}\label{thm:constrained_cond}
  For any $\beta \in (0,1]$ and $\alpha \in (0,1)$, let $k^*_\beta(\alpha) \in
  \reals_{\ge 0}$ be such that
  \begin{equation}\label{eq:N_star}
    N_\beta^{(k^*_\beta(\alpha))}(0) = \alpha.
  \end{equation}
  Such a $k^*_\beta(\alpha)$ always exists and we have the following:
  \begin{enumerate}
    \item The strategy $(f^{(k^*_\beta(\alpha))}, g^*)$ is optimal for
      Problem~\ref{prob:P2} with constraint $\alpha$. 
    \item The distortion-transmission function $D^*_\beta(\alpha)$ is
      continuous, convex and decreasing in $\alpha$ and is given by
      \begin{equation}\label{eq:D_opt}
        D^*_\beta(\alpha) = D_\beta^{(k^*_\beta(\alpha))}(0).
      \end{equation}
  \end{enumerate}
\end{theorem}

The proof of the theorem is given in Section~\ref{sec:proofs_B}.
Algorithm~\ref{algo:D_star} shows how to compute $D^*_\beta(\alpha)$.

\begin{algorithm2e}[!tb]
   \SetKwInOut{Input}{input}
   \SetKwInOut{Output}{output}
   \SetKwFor{ForAll}{forall}{do}{}
   \DontPrintSemicolon
   \Input{$\alpha \in (0,1)$, $\beta \in (0,1]$, $\varepsilon \in \reals_{>0}$}
   \Output{$D^{(k^\circ)}_\beta(\alpha)$, where $| N^{(k^\circ)}_\beta(0)- \alpha| < \varepsilon$}
   Pick $\underline k$ and $\bar k$ such that $N^{(\underline k)}_\beta(0) < \alpha < N^{(\bar k)}_\beta(0)$\;
   $k^\circ = (\underline k + \bar k)/2$\;
   \While{$| N^{(k^\circ)}_\beta(0)- \alpha| > \varepsilon$}{
     \uIf{$N^{(k^\circ)}_\beta(0) < \alpha$}{$\underline k = k^\circ$}
     \Else{ $\bar k = k^\circ$}
     $k^\circ = (\underline k + \bar k)/2$\;
   }
   \Return{$D^{(k^\circ)}_\beta(\alpha)$\;}
   \caption{Computation of $D^*_\beta(\alpha)$}\label{algo:D_star}
\end{algorithm2e}

\subsubsection{Special case of Model~B--Gauss-Markov model}\label{sec:Gauss_Markov}

In general, the optimal thresholds, and the functions $C^*_\beta(\lambda)$ and
$D^*_\beta(\alpha)$ depend on the noise distribution $\phi(\cdot)$. For the
Gauss-Markov model, the dependence on the variance $\sigma^2$ of the noise may
be quantified exactly. 

For ease of notation, we drop the dependence on $\beta$ from the notation, and
instead, show the dependence on $\sigma$. Thus, $C^*_{\sigma}(\lambda)$
denotes the optimal value for the costly communication case when the noise
variance is $\sigma^2$. Similar notation holds for other terms.

\begin{theorem}\label{thm:GM_scaling}
  For the Gauss-Markov model for Problem~\ref{prob:P1}, $k^*_{\sigma}(\lambda)
  = k^*_{1}(\lambda/a^2 \sigma^2)$ and $C^*_{\sigma}(\lambda) = \sigma^2
  C^*_{1}(\lambda/\sigma^2)$. For Problem~\ref{prob:P2}, $k^*_{\sigma}(\alpha)
  = \sigma k^*_{1}(\alpha)$ and $D^*_{\sigma}(\alpha) = \sigma^2
  D^*_{1}(\alpha)$.
\end{theorem}
The proof of the theorem is given in Section~\ref{sec:proofs_B}.

\begin{figure}[!t]
  \centering
  \begin{subfigure}[b]{0.5\textwidth}
    \includegraphics[width=0.7\linewidth]{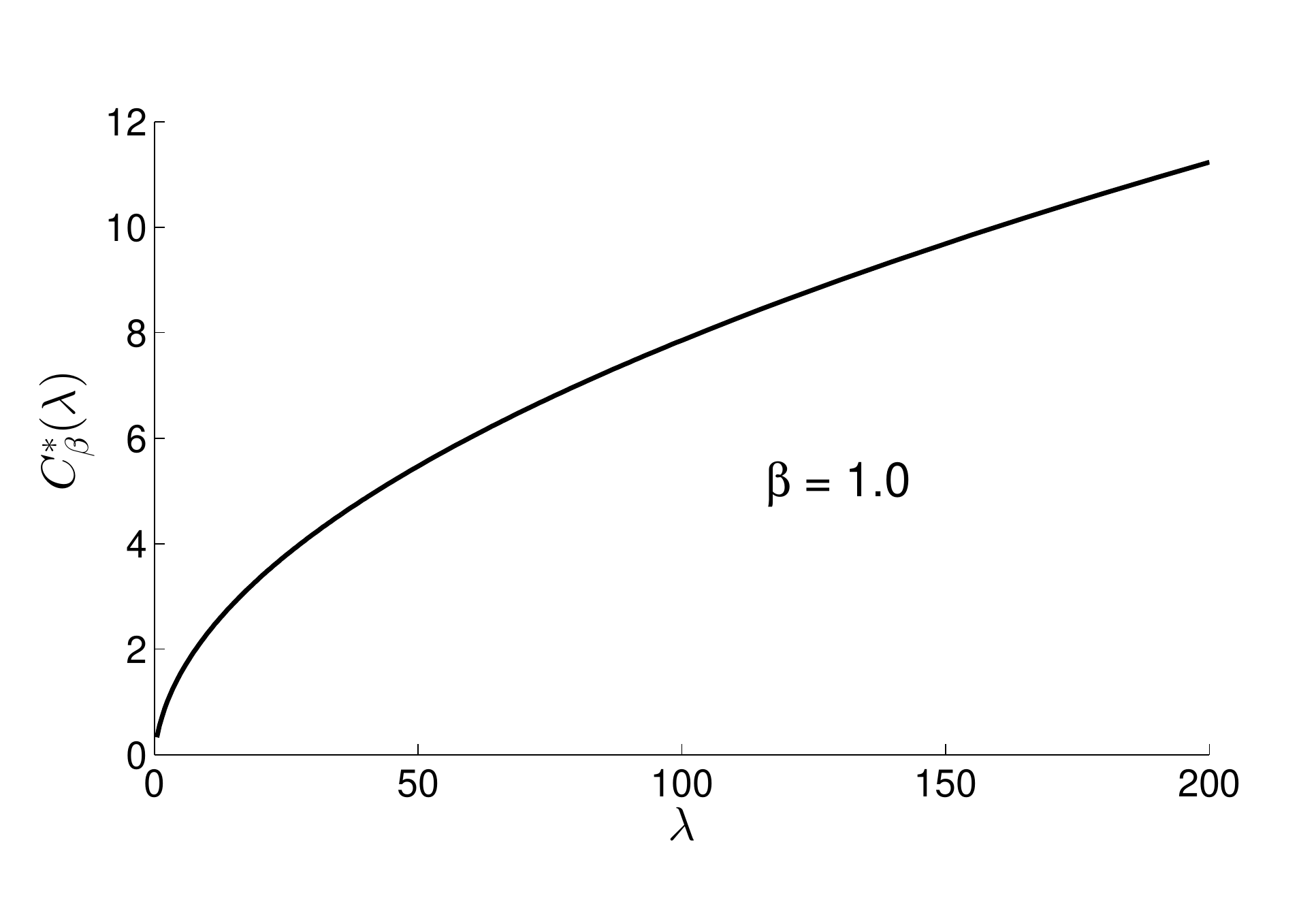}
    \caption{}
  \end{subfigure}%
  \\
  \begin{subfigure}[b]{0.5\textwidth}
    \includegraphics[width=0.7\linewidth]{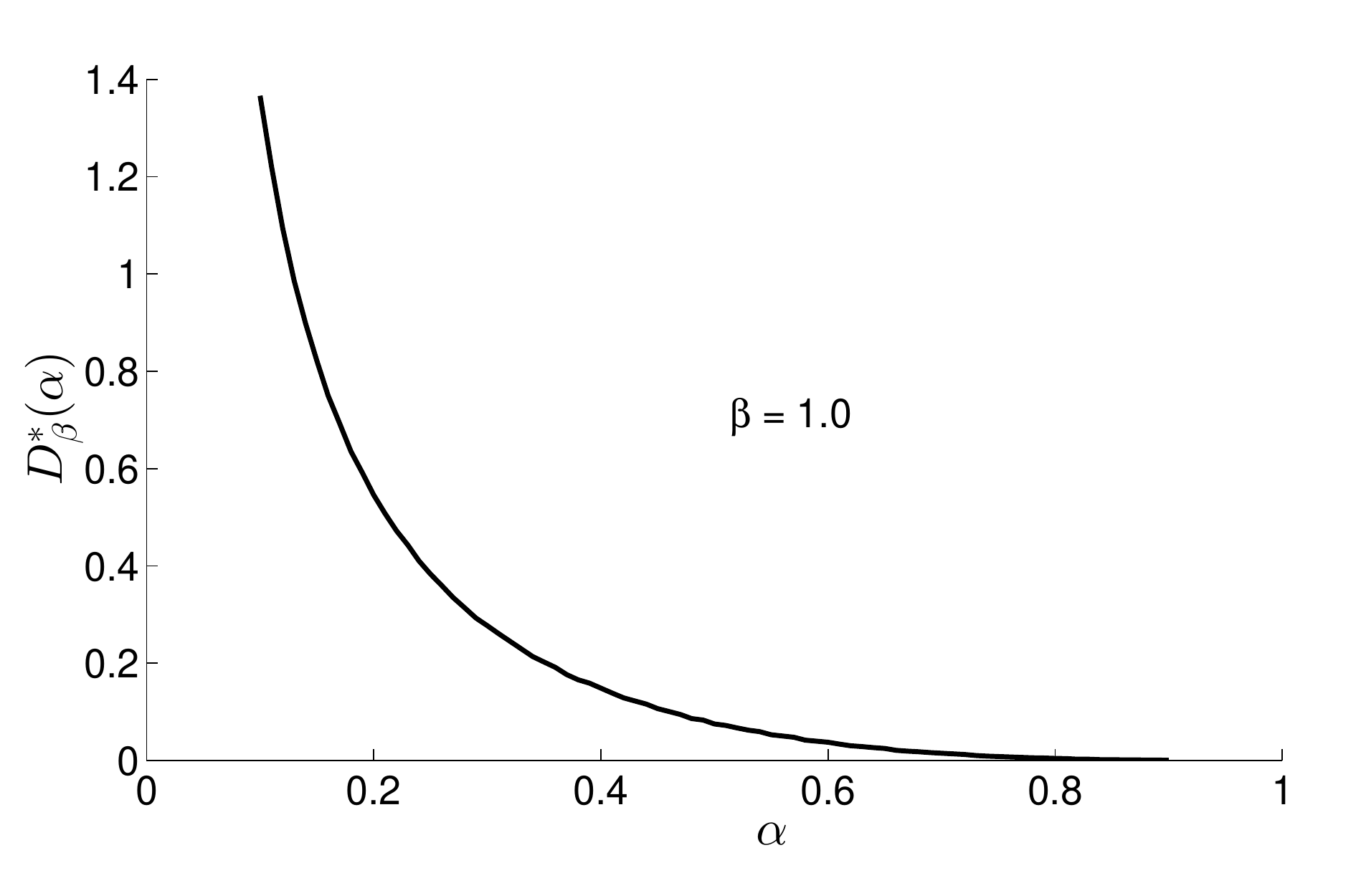}
    \caption{}
  \end{subfigure}
\caption{Gauss-Markov model ($\sigma^2 = 1$ and $a=1$): (a) optimal costly
communication cost $C^*_1(\alpha)$; (b) distortion-transmission function~$D^*_1(\alpha)$.}
\label{fig:C_star_D_star_B}
\end{figure}

An implication of the above theorem is that we only need to numerically
compute $C^*_{1}(\lambda)$ and $D^*_{1}(\alpha)$, which are shown in
Fig.~\ref{fig:C_star_D_star_B}. The optimal  total communication cost and the
distortion-transmission function for any other value $\sigma^2$ can be
obtained by simply scaling $C^*_{1}(\lambda)$ and $D^*_{1}(\alpha)$
respectively.

\subsection {An example for Model~A: symmetric birth-death Markov chain}\label{sec:example_A}

An example of a Markov process and a distortion function that satisfy Model~A is the following:

% Table to be typeset on the next page
\begin{table*}[t]
 \caption{Values of $D^{(k)}_\beta$, $N^{(k)}_\beta$ and $\lambda^{(k)}_\beta$
 for different values of $k$ and $\beta$ for the Markov chain of
 Example~\ref{ex:BD} with $p=0.3$. Note that $D^{(0)}_\beta(0) =
 D^{(1)}_\beta(0)$; therefore $\mathds K$ defined in
 Theorem~\ref{thm:DIS_AVG} equals $\integers_{> 0}$.}
 \label{tab:DNLambda}
 \begin{subtable}[h]{0.3\textwidth}
   \caption{For $\beta = 0.9$}
   \label{tab:beta09}
   \centering
   \begin{tabular}{*{4}{c}}
     \toprule
     $k$ & $D^{(k)}_\beta(0)$ & $N^{(k)}_\beta(0)$ & $\lambda^{(k)}_\beta$ \\
     \midrule
     0 & 0 & 1 & --\\
     1 & 0 & 0.5400 & 1.0989\\
     2 & 0.4576 & 0.1236 & 4.1021\\
     3 & 0.7695 & 0.0475 & 9.2839\\
     4 & 1.0066 & 0.0220 & 16.2509\\
     5 & 1.1844 & 0.0111 & 24.4478\\
     6 & 1.3130 & 0.0058 & 33.4121\\
     7 & 1.4029& 0.0031 & 42.8289\\
     8 & 1.4638 & 0.0017 & 52.5042\\
     9 & 1.5040 & 0.0009 & 62.3245\\
     10 & 1.5298 & 0.0005  & 72.2255\\
     \bottomrule
   \end{tabular}
 \end{subtable}
 \hfill
 \begin{subtable}[h]{0.3\textwidth}
   \caption{For $\beta = 0.95$}
   \label{tab:beta095}
   \centering
   \begin{tabular}{*{4}{c}}
     \toprule
     $k$ & $D^{(k)}_\beta(0)$ & $N^{(k)}_\beta(0)$ & $\lambda^{(k)}_\beta$ \\
     \midrule
     0 & 0 & 1 & --\\
     1 & 0 & 0.5700 & 1.1050\\
     2 & 0.4790 & 0.1365 & 4.3657\\
     3 & 0.8282 & 0.0565 & 10.6058\\
     4 & 1.1218 & 0.0288 & 19.9550\\
     5 & 1.3715 & 0.0163 & 32.0869\\
     6 & 1.5811 & 0.0098 & 46.4727 \\
     7 & 1.7536 & 0.0061 & 62.5651\\
     8 & 1.8927 & 0.0039 & 79.8921\\
     9 & 2.0028 & 0.0025 & 98.0854\\
     10 & 2.0884 & 0.0016 & 116.8739\\
     \bottomrule
   \end{tabular}
 \end{subtable}
 \hfill
 \begin{subtable}[h]{0.3\textwidth}
   \caption{For $\beta = 1.0$}
   \label{tab:beta1}
   \centering
   \begin{tabular}{*{4}{c}}
     \toprule
     $k$ & $D^{(k)}_\beta(0)$ & $N^{(k)}_\beta(0)$ & $\lambda^{(k)}_\beta$ \\
     \midrule
     0 & 0 & 1 & --\\
     1 & 0 & 0.6000 & 1.1111\\
     2 & 0.5000 & 0.1500 & 4.6667\\
     3 & 0.8889 & 0.0667 & 12.3810\\
     4 & 1.2500 & 0.0375 & 25.9259\\
     5 & 1.6000 & 0.0240 & 46.9697\\
     6 & 1.9444 & 0.0167 & 77.1795\\
     7 & 2.2857 & 0.0122 & 118.2222\\
     8 & 2.6250 & 0.0094 & 171.7647\\
     9 & 2.9630 & 0.0074 & 239.4737\\
     10 & 3.0000 & 0.0060 & 323.0159\\
     \bottomrule
   \end{tabular}
 \end{subtable}
\end{table*}

\begin{example}\label{ex:BD}
Consider a Markov chain of the form~\eqref{eq:state_dyn} where the pmf of $W_t$ is given by
\begin{equation*}
p_n = \begin{cases}
          p,& \text{if $|n| = 1$}\\
         1-2p, & \text{if $n=0$}\\
          0, & \text{otherwise},
            \end{cases}
\end{equation*}
where $p \in (0,\frac 13)$. The distortion function is taken as $d(e) = |e|$. 
\end{example}

\begin{figure}[!ht]
  \centering
  \includegraphics[width=\linewidth]{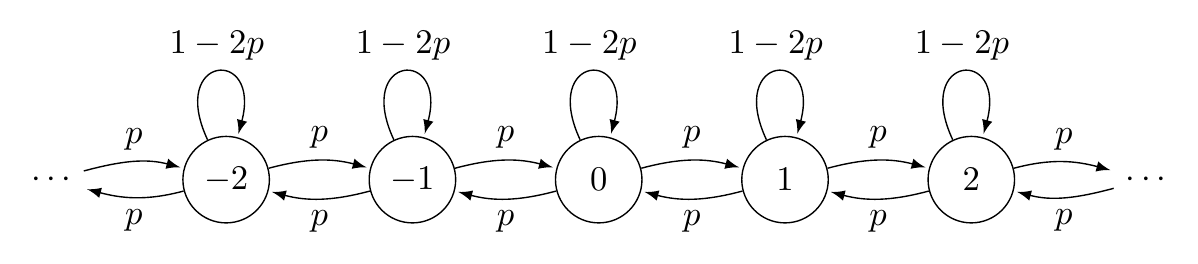}
  \caption{A birth-death Markov chain}
  \label{fig:birth-death}
\end{figure}

This Markov process corresponds to a symmetric, birth-death Markov chain defined over $\integers$ as shown in Fig.~\ref{fig:birth-death}, with the transition probability matrix is given by
\begin{equation*}
  P_{ij} = \begin{cases}
    p, & \text{if $|i-j| = 1$}; \\
    1 - 2p, & \text{if $i = j$}; \\
    0, & \text{otherwise}.
  \end{cases}
\end{equation*}

\begin{figure*}[t]
  \centering
  \begin{subfigure}[b]{0.33\textwidth}
    \includegraphics[page=1,width=\linewidth]{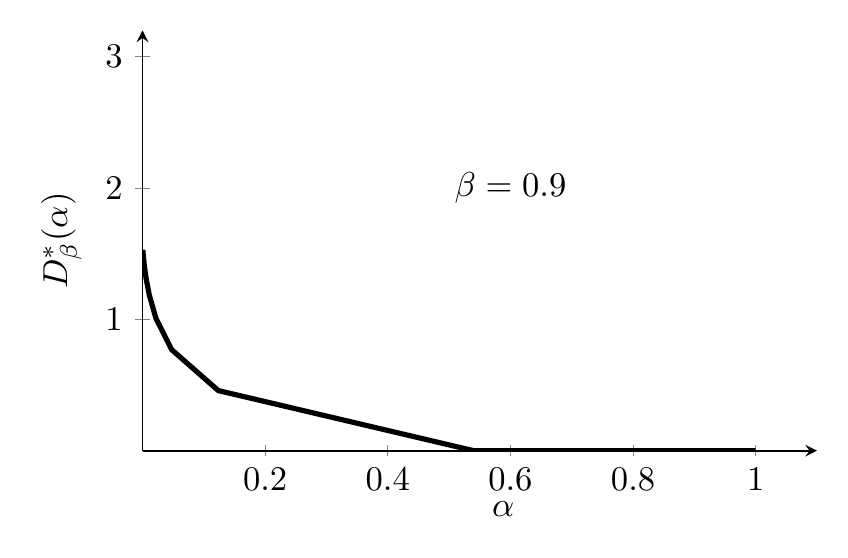}
    \caption{$D^*_\beta(\alpha)$ vs $\alpha$ for $\beta = 0.9$}
    \label{fig:D_beta09}
  \end{subfigure}%
  \hfill
  \begin{subfigure}[b]{0.33\textwidth}
    \includegraphics[page=2,width=\linewidth]{figures/chakr5}
    \caption{$D^*_\beta(\alpha)$ vs $\alpha$ for $\beta = 0.95$}
    \label{fig:D_beta095}
  \end{subfigure}%
  \hfill
  \begin{subfigure}[b]{0.33\textwidth}
    \includegraphics[page=3,width=\linewidth]{figures/chakr5}
    \caption{$D^*_\beta(\alpha)$ vs $\alpha$ for $\beta = 1.0$}
    \label{fig:D_beta1}
  \end{subfigure}
  \caption{Plots of $D^*_\beta(\alpha)$ vs $\alpha$ for different $\beta$ for
    the birth-death Markov chain of Example~\ref{ex:BD} with
  $p=0.3$.}\label{fig:DvsAlpha}
\end{figure*}

\subsubsection{Performance of a generic threshold-based strategy}
\begin{lemma}\label{lem:DN_BDMC}
  \begin{enumerate}
    \item For $\beta \in (0,1)$, 
         \begin{align*}
          D^{(k)}_\beta(0) &= \frac{\sinh(k m_\beta) - k \sinh (m_\beta) }{2 \sinh^2 (k m_\beta /2 )\sinh (m_\beta) };\\
          N^{(k)}_\beta(0) &= \frac{2 \beta p \sinh^2(m_\beta/2) \cosh(k m_\beta)}{ \sinh^2 (k m_\beta /2 )} - (1-\beta).
         \end{align*}
      \item For $\beta = 1$,
           \begin{equation*}
             D^{(k)}_1 = \frac{k^2-1}{3k};\quad
             N^{(k)}_1 = \frac{2p}{k^2};
           \end{equation*}
  \end{enumerate}
  and
  \begin{equation*}
    \lambda^{(k)}_1 = \frac {k(k+1)(k^2 + k + 1)}{6p (2k+1)}.
  \end{equation*}
\end{lemma}
The proof is given in Section~\ref{sec:proofs_example}.

\subsubsection{Optimal strategy for costly communication}

\begin{figure}[!t]
  \centering
  \includegraphics[width=0.8\linewidth]{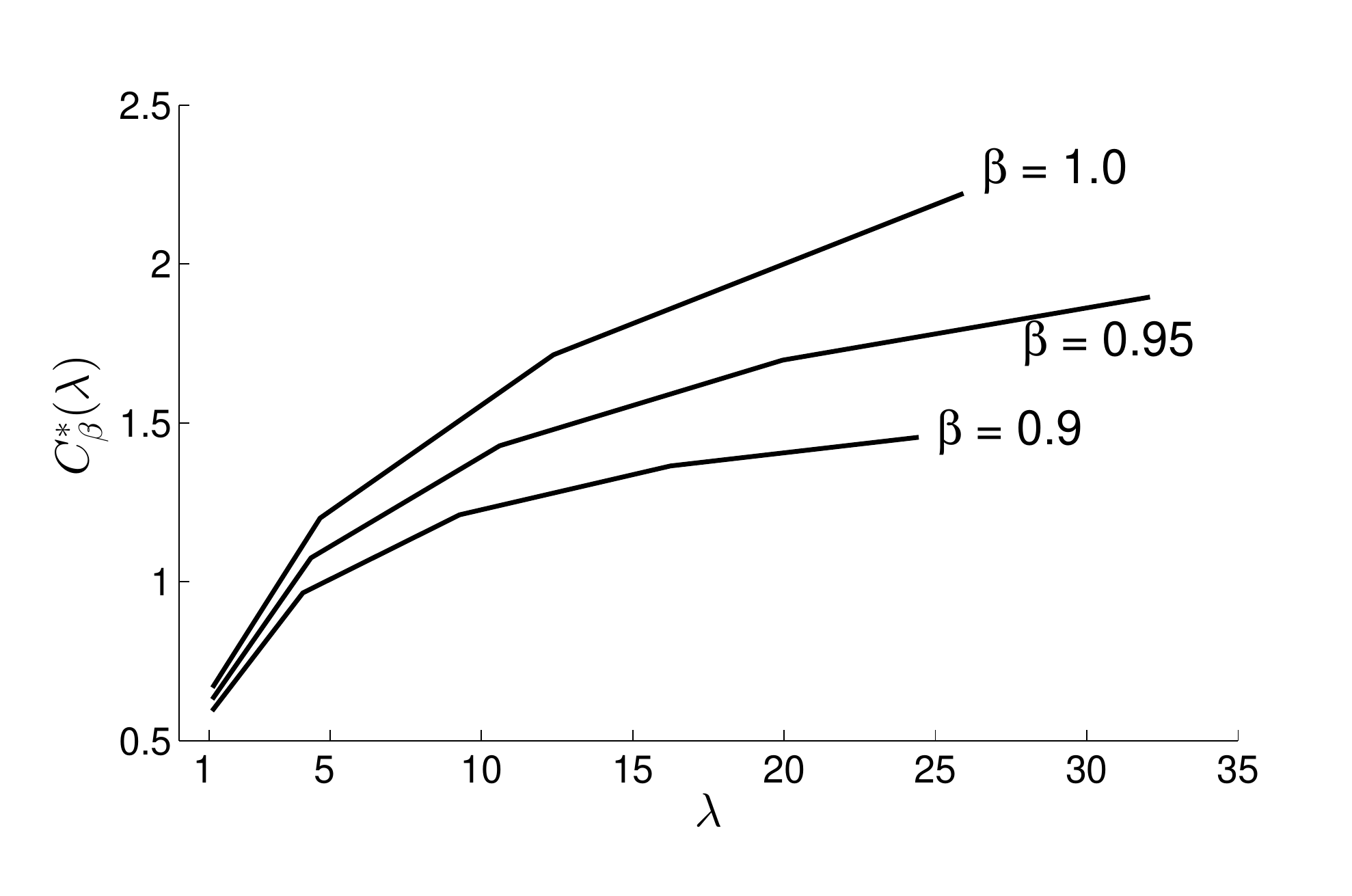}
  \caption{Plot of $C^*_\beta(\lambda)$ vs $\lambda$ for the Markov chain of
    Example~\ref{ex:BD} with $p=0.3$.}
   \label{fig:plot}
 \end{figure}

Using the above expressions for $D^{(k)}_\beta(0)$ and $N^{(k)}_\beta(0)$, we
can identify $\mathds K$ and for each $k_n \in \mathds K$, compute
$\lambda^{(k_n)}_\beta$ according to \eqref{eq:lambda}. These values
are tabulated in Table~\ref{tab:DNLambda} for different values of $\beta$ (all
for $p = 0.3$). Using Table~\ref{tab:DNLambda}, we can compute the corner
points $(\lambda^{(k_n)}_\beta, D^{(k_n)}_\beta(0) + \lambda^{(k_n)}_\beta
N^{(k_n)}_\beta(0))$ of $C^*_\beta(\lambda)$. Joining these points by straight
lines gives $C^*_\beta(\lambda)$, as shown in Fig.~\ref{fig:plot}. The optimal
strategy for a given $\lambda$ can be computed from Table~\ref{tab:DNLambda}. 

For example, for $\lambda=20$, $\beta=0.9$, we can find from
Table~\ref{tab:beta09} that $\lambda \in (\lambda^{(4)}_\beta,
\lambda^{(5)}_\beta]$. Hence, $k^*_\beta = 5$ (i.e., the strategy $f^{(5)}$ is
optimal) and the optimal total communication cost is 
\begin{equation*}
  C^*_{0.9}(20) = D^{(5)}_{0.9}(0) + 20 N^{(5)}_{0.9}(0) 
  = 1.1844+20\times 0.0111 = 1.4064.
\end{equation*}

\subsubsection{Optimal strategy for constrained communication}

Using the values in Table~\ref{tab:DNLambda}, we can also compute the corner
points $(N^{(k)}_\beta(0), D^{(k)}_\beta(0))$ of $D^*_\beta(\alpha)$. Joining
these points by straight lines gives $D^*_\beta(\alpha)$ (see
Fig.~\ref{fig:DvsAlpha}). The optimal strategy for a given $\alpha$ can be
computed from Table~\ref{tab:DNLambda}. For example, at $\alpha = 0.1$ and
$\beta = 0.9$, $k^*_\beta(\alpha)$ is the largest value of $k$ such that
$N^{(k)}_\beta(0) \ge \alpha$. Thus, from Table~\ref{tab:beta09}, we get that
$k^*=2$. Then, by~\eqref{eq:N_star},
\begin{equation*}
  \theta^* = \frac{\alpha-N^{(3)}_\beta}{N^{(2)}_\beta-N^{(3)}_\beta} = 0.6899.
\end{equation*}
Let $f^* = (f^{(2)}, f^{(3)}, \theta^*)$. Then the Bernoulli randomized simple strategy $(f^*, g^*)$ is optimal for Problem~\ref{prob:P2} for $\beta \in (0,1)$. Furthermore, by~(\ref{eq:D_opt_rand}), $D^*_\beta(\alpha) = 0.5543$. 

\section{Salient features and discussion}\label{sec:salient_fea_discussion}
\subsection{Comparison with periodic and randomized strategies}
\label{subsec:per_rand}

In our model, we assume that the transmission decision depends on the state
of the Markov process. In some of the remote estimation literature, it is
assumed that the transmission schedule does not depend on the state of the
Markov process. Two such commonly used strategies are:
\begin{enumerate}
  \item Periodic transmission strategy with period $T$:
    \begin{equation*}
      U_t = f_p(t \text{ mod } T),
    \end{equation*}
    where $\sum_{t=0}^{T-1} f_p(t) = 1/\alpha$.
  \item Random transmission strategy: 
    \begin{equation*}
      U_t = \begin{cases}
        1, & \mbox{w.p. $\alpha$}\\
        0, & \mbox{w.p. $1-\alpha$}.
      \end{cases}
    \end{equation*}
\end{enumerate}
Below, we compare the performance of the threshold-based strategy with these
two strategies for the for the long-term average setup for
Problem~\ref{prob:P2} for Model~B with $a=1$.

\subsubsection{Performance of the periodic strategy}
\label{subsubsec:perform_per}

In general, the performance of a periodic transmission strategy depends on
the choice of transmission function $f_p$. For ease of calculation we
consider the values of $(\alpha,T)$ for which $f_p$ is unique.

\begin{enumerate}
  \item $\alpha = 1/T$, $T \in \integers_{>0}$, i.e., the transmitter remains
    silent for $(T-1)$ steps and then transmits once. The expected distortion
    in this case is 
    \begin{align*}
      D_{\text{per}}(\alpha) 
      &= \frac{1}{T} \EXP \Big[ \sum_{t=0}^{T-1} E^2_t \Big]\\
      &\stackrel{(a)} = \frac{1}{T} \EXP \Big[ \sum_{t=0}^{T-1} t \sigma^2 \Big] 
      = \frac{1}{T} \frac{(T-1)T}{2} \sigma^2 = \frac{\sigma^2}{2} \big(\frac{1}{\alpha} -1\big),
    \end{align*}
    where $(a)$ uses $E_t = W_0+ W_1+ W_2+ \cdots+ W_{t-1}$.
  \item $\alpha = (T-1)/T$, $T \in \integers_{> 0}$, i.e., the transmitter
    remains silent for 1 step and then transmits for $(T-1)$ steps. The
    expected distortion in this case is
    \begin{equation*}
      D_{\text{per}}(\alpha) = \frac{1}{T} \EXP [E^2_1] 
      = \frac{\sigma^2}{T} = \sigma^2(1-\alpha).
    \end{equation*}
\end{enumerate}

\subsubsection{Performance generic stationary transmission strategy}

Next, we derive an expression of $D_\beta(f, g^*)$ for arbitrary stationary
transmission strategy $f$ (that does not use the value of the state $E_t$ to
determine when to transmit; so the receiver is the same as in
Theorem~\ref{thm:inf_hor}) for the long-term average setup for Model~B when
$a=1$.

\begin{proposition}\label{prop:D_arbit_stationary}
  For $\beta = 1$ and $a=1$ in Model~B, let $f$ be an arbitrary stationary
  transmission strategy. Let $\tau$ denote the stopping time of the first
  transmission under $f$. Then
  \begin{equation*}
    D_1(f, g^*) 
    = \frac{\sigma^2}{2} \Big[\frac{\EXP (\tau^2)}{\EXP (\tau)} - 1\Big].
  \end{equation*}
\end{proposition}

\begin{IEEEproof}
  For any $t < \tau$, $E_t = W_0^2 + \dots + W_{t-1}^2$. Therefore,
  $\EXP[E^2_t] = t \sigma^2$ and define $\hat L(t) = \sum_{s=1}^{t -1 }
  \EXP[E_s^2] = \frac 12 t(t-1) \sigma^2$. Now, $L_1(0) = \EXP [\hat{L}
  (\tau)] = (\sigma^2/2) [\EXP(\tau^2) - \EXP(\tau)]$ and $M_1(0) =
  \EXP(\tau)$. By using the same argument as in the proof of
  Theorem~\ref{thm:renewal_DNC}, we get $D_1(f, g^*) = L_1(0)/M_1(0)$,
  which implies the result.
\end{IEEEproof}

\subsubsection{Performance of randomized transmission strategy}

For the randomized strategy defined above, $\tau$ is a
$\text{Geom}_1(\alpha)$ random variable. Therefore, $\EXP (\tau^2) =
2/\alpha^2 - 1/\alpha$ and $\EXP (\tau) = 1/\alpha$. Hence, following
Proposition~\ref{prop:D_arbit_stationary}, we have 
\begin{equation*}
  D_{\text{rand}}(\alpha) = \sigma^2 \big[\frac{1}{\alpha} -1 \big].
\end{equation*}

\begin{figure}[!ht]
  \centering
  \includegraphics[width=0.8\linewidth]{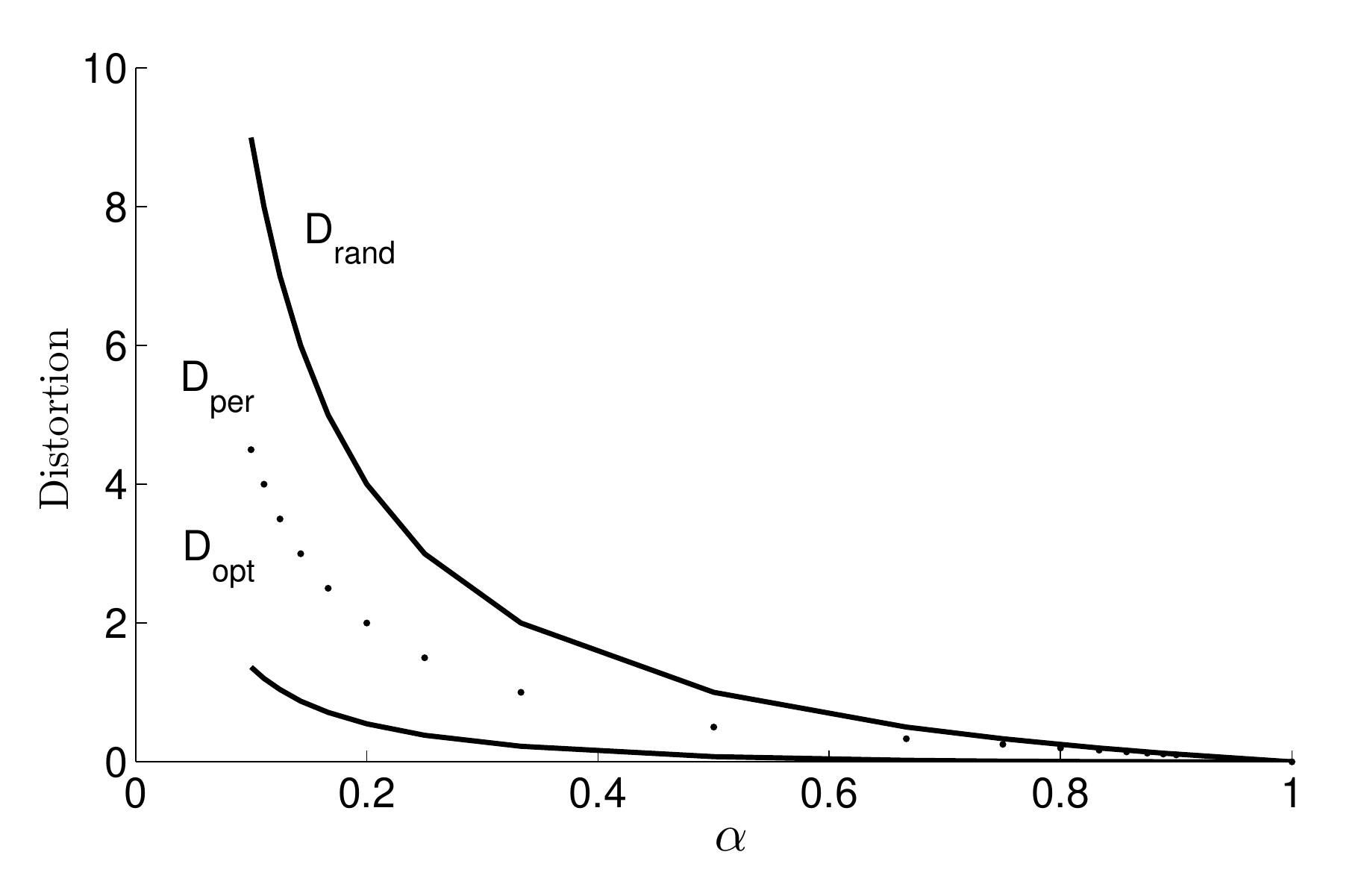}
  \caption{Comparison of the performances of the threshold-based startegy
    (denoted by $D_{\text{opt}}$) with periodic and randomized transmission
    strategies (denoted by $D_{\text{per}}$ and $D_{\text{rand}}$,
  respectively) for a Gauss-Markov process with $a = 1$ and $\sigma^2 = 1$.}
  \label{fig:D_th_per_rand}
\end{figure}

Fig.~\ref{fig:D_th_per_rand} shows that threshold-based startegy performs
considerably well compared to the periodic transmission strategy and the
randomized transmission strategy. 

\subsection{Discussion on deterministic implementation}

The optimal strategy shown in Theorem~\ref{thm:constrained_f_star_D} chooses
a randomized action in states $\{-k^*,k^*\}$. It is also possible to identify
deterministic (non-randomized) but time-varying strategies that achieve the
same performance. We describe two such strategies for the long-term average
setup.

\subsubsection{Steering strategies}

Let $a^0_t$ (respectively, $a^1_t$) denote the number of times the action
$u_t=0$ (respectively, the action $u_t=1$) has been chosen in states
$\{-k^*,k^*\}$ in the past, i.e.
\begin{equation*}
  a^i_t = \sum_{s=0}^{t-1} \IND \{|E_s| = k^*,\,u_s=i\}, \quad i \in \{0,1\}.
\end{equation*}
Thus, the empirical frequency of choosing action $u_t=i$, $i \in \{0,1\}$, in
states $\{-k^*,k^*\}$ is $a^i_t/(a^0_t+a^1_t)$. A steering strategy compares
these empirical frequencies with the desired randomization probabilities
$\theta^0=1-\theta^*$ and $\theta^1 = \theta^*$ and chooses an action that
\emph{steers} the empirical frequency closer to the desired randomization
probability. More formally, at states $\{-k^*,k^*\}$, the steering
transmission strategy chooses the action
\begin{equation*}
  \arg\min_i \Big\{\theta^i - \frac{a^i_t+1}{a^0_t+a^1_t+1}\Big\}
\end{equation*}
in states $\{-k^*,k^*\}$ and chooses deterministic actions according to $f^*$
(given in (\ref{eq:f_star})) in states except $\{-k^*,k^*\}$. Note that the
above strategy is deterministic (non-randomized) but depends on the history
of visits to states $\{-k^*,k^*\}$. Such strategies were proposed in
\cite{Feinberg}, where it was shown that the steering strategy descibed above
achieves the same performance as the randomized startegy $f^*$ and hence is
optimal for Problem~\ref{prob:P2} for $\beta=1$. Variations of such steering
strategies have been proposed in \cite{shwartz1986optimal,ma1990stochastic},
where the adaptation was done by comparing the sample path average cost with
the expected value (rather than by comparing empirical frequencies).

\subsubsection{Time-sharing strategies}

Define a cycle to be the period of time between consecutive visits of process
$\{E_t\}_{t=0}^\infty$ to state zero. A time-sharing strategy is defined by a
series $\{(a_m,b_m)\}_{m=0}^\infty$ and uses startegy $f^{(k^*)}$ for the
first $a_0$ cycles, uses startegy $f^{(k^*+1)}$ for the next $b_0$ cycles,
and continues to alternate between using startegy $f^{(k^*)}$ for $a_m$
cycles and strategy $f^{(k^*+1)}$ for $b_m$ cycles. In particular, if
$(a_m,b_m) = (a,b)$ for all $m$, then the time-sharing strategy is a periodic
strategy that uses $f^{(k^*)}$ $a$ cycles and $f^{(k^*+1)}$ for $b$ cycles. 

The performance of such time-sharing strategies was evaluated in
\cite{altman_shwartz}, where it was shown that if the cycle-lengths of the
time-sharing strategy are chosen such that,
\begin{align*}
  \lim_{M \to \infty} \frac{\sum_{m=0}^M a_m}{\sum_{m=0}^M (a_m+b_m)} &= \frac{\theta^* N^{(k^*)}_1}{\theta^* N^{(k^*)}_1 + (1-\theta^*) N^{(k^*+1)}_1}\\
  & = \frac{\theta^* N^{(k^*)}_1}{\alpha},
\end{align*}
then the time-sharing strategy $\{(a_m,b_m)\}_{m=0}^\infty$ achieves the same
performance as the randomized strategy $f^*$ and hence, is optimal for
Problem~\ref{prob:P2} for $\beta =1$.

\section{Proof of the structural result: Theorem~\ref{thm:inf_hor}}
\label{sec:proof_struct}
\subsection{Finite horizon setup}\label{subsec:finite_hor}
A finite horizon version of Problem~\ref{prob:P1} has been investigated
in~\cite{NayyarBasarTeneketzisVeeravalli:2013} (for Model A) and
in~\cite{LipsaMartins:2011,MH2012} (for Model~B), where the structure of the
optimal transmission and estimation strategy was established. 

\begin{theorem}\label{thm:fin_hor}
  \cite{LipsaMartins:2011,NayyarBasarTeneketzisVeeravalli:2013,MH2012}
  For both Models~A and~B, for a finite horizon version of
  Problem~\ref{prob:P1}, we have the following.
  \begin{enumerate}
    \item \emph{Structure of optimal estimation strategy}: the estimation
      strategy defined in Theorem~\ref{thm:inf_hor} is optimal.
    \item \emph{Structure of optimal transmission strategy}: define $E_t$ as
      in Theorem~\ref{thm:inf_hor}. Then there exist threholds
      $\{k_t\}_{t=1}^T$ such that the transmission strategy
      \begin{equation} \label{eq:transmit-opt}
        U_t \DEFINED f_t(E_t) = \begin{cases}
          1, & \text{if $|E_t| \ge k_t$}; \\
          0, & \text{if $|E_t| < k_t$} 
        \end{cases}
      \end{equation}
      is optimal. 
  \end{enumerate}
\end{theorem}

The above structural results were obtained in~\cite[Theorems~2
and~3]{NayyarBasarTeneketzisVeeravalli:2013} for Model A and in~\cite[Theorem
1]{LipsaMartins:2011} and~\cite[Lemmas 1, 3 and 4]{MH2012} of Model B.

\begin{remark}\label{rem:countable_support}
  The results in~\cite{NayyarBasarTeneketzisVeeravalli:2013} were derived
  under the assumption that $\{W_t\}$ has finite support. These results can
  be generalized for $\{W_t\}$ having countable support using ideas
  from~\cite{WangWooMadiman:2014}. For that reason, we state
  Theorem~\ref{thm:fin_hor} without any restriction on the support of
  $\{W_t\}$. 
  See the supplementary document for the generalization of~\cite[Theorems~2 and~3]{NayyarBasarTeneketzisVeeravalli:2013} to $\{W_t\}$ with countable support.
\end{remark}

\subsection{Infinite horizon setup}\label{subsec:inf_hor}
In a general real-time communication system, the optimal estimation strategy
depends on the choice of the transmission strategy and vice-versa.
Theorem~\ref{thm:fin_hor} shows that when the noise process and the
distortion function satisfy appropriate symmetry assumptions, the optimal
estimation strategy can be specified in closed form. Consequently, we can fix
the estimation strategy to be of the above form and consider the optimization
problem of identifying the best transmission strategy. This optimization
problem has a single decision maker---the transmitter---and we use techniques
from centralized stochastic control to solve it. Since the optimal estimation
strategy is time-homogeneous, one expects the optimal transmission strategy
(i.e., the choice of the optimal thresholds $\{k_t\}_{t=0}^\infty$) to be
time-homogeneous as well. The technical difficulty in establishing such a
result is that the state space is not compact and the distortion function may
be unbounded. 

To prove Theorem~\ref{thm:inf_hor}, we proceed as follows:
\begin{enumerate}
  \item We show that the result of the theorem is true for $\beta \in (0,1)$
    and the optimal strategy is given by an appropriate dynamic program. 
  \item We show that for the discounted setup, the value function of the
    dynamic program is even and increasing on $\mathds X$.
  \item For $\beta =1$, we use the vanishing discount approach to show that
    the optimal strategy for the long-term average cost setup may be
    determined as a limit to the optimal strategy for the discounted cost
    setup is the discount factor $\beta \uparrow 1$.
\end{enumerate}

\subsubsection{The discounted setup}

\begin{lemma}\label{lem:DP_A}
  In Model~A. an optimal transmission strategy is given by the unique and
  bounded solution of the following dynamic program: for all $e \in
  \integers$,
  \begin{align}
    V_\beta(e;\lambda) &= \min \Big[(1-\beta)\lambda + \beta \sum_{w \in
    \integers} p_{w} V_\beta(w;\lambda),
    \notag \\
    &\quad (1-\beta)d(e) + 
    \beta \sum_{w \in \integers} p_{w} V_\beta(ae+w;\lambda)\Big].
    \label{eq:beta_optimality_A}
    \end{align}
\end{lemma}

\begin{IEEEproof}
  When $d(\cdot)$ is bounded, the per-step cost $c(e,u) \DEFINED
  (1-\beta)[\lambda u + d(e)(1-u)]$, $u \in \{0,1\}$, for a given $\lambda$
  is also bounded and hence according to \cite[Proposition 4.7.1, Theorem
  4.6.3]{Sennott:book}, there exists the unique and bounded solution
  $V_\beta(e;\lambda)$ of the dynamic program~\eqref{eq:beta_optimality_A}.  

  When $d(\cdot)$ is unbounded, then for any communication cost $\lambda$, we
  first define $e_0 \in \integers_{\ge 0} < \infty$ as:
  \begin{equation*}
    e_0 \DEFINED \min \Big\{e \,:\, d(e) \ge \frac{\lambda}{1-\beta} \Big\}.
  \end{equation*}

  Now, for any state $e$, $|e| > e_0$, the per-step cost $(1-\beta)d(e)$ of
  not transmitting is greater then the cost of transmitting at each step in
  the future, which is given by $(1-\beta) \sum_{t=0}^\infty \beta^t \lambda
  = \lambda$. Thus, the optimal action is to transmit, i.e., $f^* (e)= 1$.
  Hence, the dynamic program can be written as
  \begin{equation*}
    V_\beta(e; \lambda) = \min \{V^0_\beta(e;\lambda), V^1_\beta(e;\lambda)\},
  \end{equation*}
  where
  \begin{align*}
    V^0_\beta(e;\lambda) &= (1-\beta)d(e) + \beta \sum_{w \in \integers} p_{w} V_\beta(ae+w;\lambda),\\
    V^1_\beta(e;\lambda) &= (1-\beta)\lambda + \beta \sum_{w \in \integers} p_{w} V_\beta(w;\lambda).
  \end{align*}

  Let $\mathcal{E}^* \DEFINED \{e \, :\, |e|\ge e_0\}$. Then, for all $e \in
  \mathcal{E}^*$, $V_\beta(e;\lambda)$ is constant. Thus,
  \eqref{eq:beta_optimality_A} is equivalent to a finite-state Markov
  decision process with state space $\{-e_0+1,\cdots, e_0-1\}\cup e^*$ (where
  $e^*$ is a generic state for all states in the set $\mathcal{E}^*$). Since
  the state space is now finite, the dynamic
  program~\eqref{eq:beta_optimality_A} has a unique and bounded
  time-homogeneous solution by the argument given for bounded $d(\cdot)$. 
\end{IEEEproof}
 
\begin{lemma}\label{lem:DP_B}
  In Model~B, an optimal transmission strategy is given by the unique and
  bounded solution of the following dynamic program: for all $e \in \reals$,
  \begin{multline}
    V_\beta(e;\lambda) =  
    \min \Big[ (1-\beta) \lambda +  
      \beta \int_\reals \phi(w) V_\beta(w;\lambda)dw, \\
      \qquad (1-\beta) d(e) +  
      \beta \int_\reals \phi(w) V_\beta(ae+w;\lambda)dw
      \Big].
    \label{eq:beta_optimality_B}  
  \end{multline}
\end{lemma}
\begin{IEEEproof}
  When $d(\cdot)$ is bounded, the per-step cost $c(e,u)$, as defined in part
  (a), for a given $\lambda$ is also bounded. Let $K = (1-\beta) \sup_{e \in
  \reals} \{  d(e) \}$. Then, the strategy `always transmit'
  satisfies~\cite[Assumption 4.2.2]{LermaLasserre} with $V_\beta(e;\lambda)
  \le K/(1-\beta)$. Also, $\lambda$, $d(\cdot)$ and  $\phi(\cdot)$
  satisfy~\cite[Assumption 4.2.1]{LermaLasserre}. Hence, the above dynamic
  program has a unique and bounded solution due to~\cite[Theorem
  4.2.3]{LermaLasserre}.

  When $d(\cdot)$ is unbounded, define $e_0$ and $e^*$ as in the proof of
  Lemma~\ref{lem:DP_A}. By an argument similar to that in the proof of
  Lemma~\ref{lem:DP_A}, we can restrict the state space
  of~\eqref{eq:beta_optimality_B}   to $[-e_0,e_0]\cup e^*$. Hence, the state
  space is compact and on this state space $d(\cdot)$ is bounded. Thus, the
  dynamic program~\eqref{eq:beta_optimality_B} has a unique and bounded
  solution by the argument given for bounded $d(\cdot)$.
\end{IEEEproof}

\begin{IEEEproof}[Proof of Theorem~\ref{thm:inf_hor} for $\beta \in (0,1)$]
  The structure of the optimal strategies follows from
  Theorem~\ref{thm:fin_hor}. The optimal thresholds are time invariant because
  the corresponding dynamic programs \eqref{eq:beta_optimality_A} and
  \eqref{eq:beta_optimality_B} have a unique fixed point.
\end{IEEEproof}

\subsubsection{Properties of the value function}
\begin{proposition}\label{prop:V_plus_minus}
  For any $a \in \mathds X_{> 0}$, consider the two Markov processes
  $\{X^{(+)}_t\}_{t=0}^\infty$ and $\{X^{(-)}_t\}_{t=0}^\infty$ such that
  $X^{(+)}_0 = X^{(-)}_0 = 0$ and 
  \begin{equation*}
    X^{(+)}_{t+1} = a X^{(+)}_t + W_t
    \quad  \text{and} \quad 
    X^{(-)}_{t+1} = -a X^{(-)}_t + W_t.
  \end{equation*}
  Let $V^{(+)}_\beta$ and $V^{(-)}_\beta$ be the value functions
  corresponding to $\{X^{(+)}_t\}_{t=0}^\infty$ and
  $\{X^{(+)}_t\}_{t=0}^\infty$. Then
  \begin{equation*}
    V^{(+)}_\beta (e) = V^{(-)}_\beta (e), \quad \forall e.
  \end{equation*}
  Therefore, if $k$ is an optimal threshold for $\{X^{(+)}_t\}_{t=0}^\infty$
  then $k$ is also optimal for $\{X^{(-)}_t\}_{t=0}^\infty$.
\end{proposition}
See Appendix~\ref{app:V_plus_minus_EI} for the proof.

\begin{remark}\label{rem:a_pos} 
  As a consequence of the above proposition, we can restrict attention to $a >
  0$ while proving the properties of the value function $V_\beta(\cdot)$.
\end{remark}

\begin{proposition} \label{prop:EI}
  For any $\lambda > 0$ and $\beta \in (0,1)$, the value functions
  $V_\beta(\cdot;\lambda)$ given by~\eqref{eq:beta_optimality_A}
  and~\eqref{eq:beta_optimality_B} are even and increasing on $\mathds X_{\ge
  0}$.
\end{proposition}

See Appendix~\ref{app:V_plus_minus_EI} for the proof.

\subsubsection{The long-term average setup}

\begin{proposition}
  \label{prop:SEN}
  For any $\lambda \ge 0$, the value function $V_\beta(\cdot; \lambda)$ for
  Models~A and~B, as given by~\eqref{eq:beta_optimality_A}
  and~\eqref{eq:beta_optimality_B} respectively, satisfy the following SEN
  conditions of~\cite{Sennott:book,LermaLasserre}:
   \begin{enumerate}
    \item[(S1)] There exists a reference state $e_0 \in \mathds X$ and a
      non-negative scalar $M_\lambda$ such that $V_\beta(e_0,\lambda) <
      M_\lambda$ for all $\beta \in (0,1)$.
    \item[(S2)] Define $h_\beta(e;\lambda) = (1-\beta)^{-1}[
        V_\beta(e;\lambda) - V_\beta(e_0;\lambda)]$. There exists a function
        $K_\lambda : \integers  \to \reals$ such that $h_\beta(e;\lambda) \le
        K_\lambda(e)$ for all $e \in \mathds X$ and $\beta \in (0,1)$. 
    \item[(S3)] There exists a non-negative (finite) constant $L_\lambda$
      such that $-L_\lambda \le h_\beta(e;\lambda)$ for all $e \in \mathds X$
      and $\beta \in (0,1)$. 
  \end{enumerate}
  Therefore, if $f_\beta$ denotes an optimal strategy for $\beta \in (0,1)$,
  and $f_1$ is any limit point of $\{f_\beta\}$, then $f_1$ is optimal for
  $\beta =1$.
\end{proposition}

\begin{IEEEproof}
  Let $V^{(0)}_\beta(e,\lambda)$ denote the value function of the `always
  transmit' strategy. Since $V_\beta(e,\lambda) \le V^{(0)}_\beta(e,\lambda)$
  and $V^{(0)}_\beta(e,\lambda) = \lambda$, (S1) is satisfied with $M_\lambda
  = \lambda$.

  We show (S2) for Model~B, but a similar argument works for Model~A as well.
  Since not transmitting is optimal at state 0, we have
  \begin{equation*}
    V_\beta(0,\lambda) = \beta \int_{-\infty}^\infty \phi(w) V_\beta(w,\lambda)dw.
  \end{equation*}
  Let $V^{(1)}_\beta(e,\lambda)$ denote the value function of the strategy
  that transmits at time 0 and follows the optimal strategy from then on. Then
  \begin{align}
    V^{(1)}_\beta(e,\lambda) &= (1-\beta) \lambda + \beta \int_{-\infty}^\infty \phi(w) V_\beta(w,\lambda)dw
    \notag \\
    &= (1-\beta)\lambda + \beta V_\beta(0,\lambda)
    \label{eq:proofVB}
  \end{align}
  Since $V_\beta(e,\lambda) \le V^{(1)}_\beta(e,\lambda)$ and
  $V_\beta(0,\lambda) \ge 0$, from~\eqref{eq:proofVB} we get that
  $(1-\beta)^{-1} [V_\beta(e,\lambda) - V_\beta(0,\lambda)] \le \lambda$.
  Hence (S2) is satisfied with $K_\lambda(e) = \lambda$.

  By Proposition~\ref{prop:EI}, $V_\beta(e,\lambda) \ge V_\beta(0,\lambda)$,
  hence (S3) is satisfied with $L_\lambda = 0$.
\end{IEEEproof}

\begin{IEEEproof}[Proof of Theorem~\ref{thm:inf_hor} for $\beta = 1$]
  Since the value function $V_\beta(\cdot,\lambda)$ satisfies the SEN
  conditions for reference state $e_0=0$, the optimaity of the threshold
  strategy for long-term average setup follows from~\cite[Theorem
  7.2.3]{Sennott:book} for Model~A and~\cite[Theorem 5.4.3]{LermaLasserre}
  for Model~B, respectively.
\end{IEEEproof}

\section{Proof of Theorem~\ref{thm:renewal_DNC}}\label{sec:proof_thm_renewal}

\subsection{Preliminary results}
Define operator $\mathcal{B}$ as follows:
\begin{itemize}
  \item \textbf{Model A}: For any $v: \integers \to \reals$, define operator
    $\mathcal{B}$ as
    \begin{equation*}
      [\mathcal{B}v](e) \DEFINED \sum_{w=-\infty}^\infty p_w v(ae+w), \quad \forall e \in \integers.
    \end{equation*}
    Or, equivalently,
    \begin{equation*}
      [\mathcal{B}v](e) \DEFINED \sum_{n=-\infty}^\infty p_{n-ae} v(n), \quad \forall e \in \integers.
    \end{equation*}
  \item \textbf{Model B}: For any bounded $v: \reals \to \reals$, define
    operator $\mathcal{B}$ as
    \begin{equation*}
      [\mathcal{B}v](e) \DEFINED \int_\reals \phi(w) v(ae+w)dw, \quad \forall e \in \reals.
    \end{equation*}
    Or, equivalently,
    \begin{equation*}
      [\mathcal{B}v](e) \DEFINED \int_\reals \phi(n-ae) v(n)dn, \quad \forall e \in \reals.
    \end{equation*}
\end{itemize}

As discussed in Remark~\ref{rem:reset}, the error process $\{E_t\}_{t=0}^\infty$ is a controlled Markov process. Therefore, the functions $D^{(k)}_\beta$ and $N^{(k)}_\beta$ may be thought as value functions when strategy $f^{(k)}$ is used. Thus, they satisfy the following fixed point equations: for $\beta \in (0,1)$,
\begin{align}\label{eq:D_fixed_pt}
  D^{(k)}_\beta(e) &=\begin{cases}
    \beta [\mathcal{B} D^{(k)}_\beta](0), & \mbox{if $|e| \ge k$}\\
    (1-\beta)d(e) + \beta [\mathcal{B} D^{(k)}_\beta](e), & \mbox{if $|e| < k$},
  \end{cases}\\ \label{eq:N_fixed_pt}
  N^{(k)}_\beta(e) &=\begin{cases}
    (1-\beta) + \beta[\mathcal{B} N^{(k)}_\beta](0), & \mbox{if $|e| \ge k$}\\
    \beta [\mathcal{B} N^{(k)}_\beta](e), & \mbox{if $|e| < k$}.
  \end{cases}
\end{align}

\begin{lemma}\label{lem:DN_lim}
  For $\beta \in (0,1]$, \eqref{eq:D_fixed_pt} and~\eqref{eq:N_fixed_pt} have
  unique and bounded solutions $D^{(k)}_\beta(e)$ and $N^{(k)}_\beta(e)$ that 
  \begin{enumerate}
    \item are even and increasing (on $\mathds X_{\ge 0}$) in $e$ for all~$k$,
    \item satisfy the SEN conditions (see Proposition~\ref{prop:SEN}) and
      therefore 
      \begin{equation*}
        D^{(k)}_1(e) = \lim_{\beta \uparrow 1} D^{(k)}_\beta (e)
        \quad \text{and}\quad
        N^{(k)}_1(e) = \lim_{\beta \uparrow 1} N^{(k)}_\beta (e).
      \end{equation*}
    \item $D^{(k)}_\beta(e)$ is increasing in $k$ for all~$e$ and
      $N^{(k)}_\beta(e)$ is strictly decreasing in $k$ for all~$e$. 
  \end{enumerate}
\end{lemma}

The proofs of 1) and 2) follow from the arguments similar to those of Section~\ref{sec:proof_struct} and are therefore omitted. The proof of 3) is given in~Appendix~\ref{app:D_inc}.

\subsection{Proof of Theorem~\ref{thm:renewal_DNC}}

We prove the result for the discounted cost setup, $\beta \in (0,1)$. The
result extends to the long-term average cost setup, $\beta =1$, by using the
vanishing discount approach similar to the argument given in
Section~\ref{sec:proof_struct}.

We first consider the case $k=0$. In this case, the recursive definition of
$D^{(k)}_\beta$ and $N^{(k)}_\beta$, given by~\eqref{eq:D_fixed_pt}
and~\eqref{eq:N_fixed_pt}, simplify to the following:
\begin{equation*}
  D^{(0)}_\beta(e) = \beta [\mathcal{B}D^{(0)}_\beta](0);
\end{equation*}
and
\begin{equation*}
  N^{(0)}_\beta(e) = (1-\beta) + \beta [\mathcal{B}N^{(0)}_\beta](0).
\end{equation*}

It can be easily verified that $D^{(0)}_\beta(e) = 0$ and $N^{(0)}_\beta(e) =
1$, $e \in \mathds X$, satisfy the above equations. Also,
$C^{(0)}_\beta(e;\lambda) = C_\beta(f^{(0)},g^*;\lambda) = \lambda$. This
proves the first part of the proposition.

For $k>0$, let $\tau^{(k)}$ denote the stopping time when the Markov process
in both Model A and B starting at state $0$ at time $t=0$ leaves the set
$S^{(k)}$. Note that $\tau^{(0)} = 1$ and $\tau^{(\infty)} = \infty$.

Then,
\begin{align}
  L^{(k)}_\beta(0) &= \EXP \Big[\sum_{t=0}^{\tau^{(k)}-1} \beta^t d(E_t)\,\Bigm|\, E_0=0 \Big]
  \label{eq:proof_L_k} \\
  M^{(k)}_\beta(0) &= \EXP \Big[\sum_{t=0}^{\tau^{(k)}-1} \beta^t \,\Bigm|\, E_0=0 \Big] = \frac{1-\EXP [\beta^{\tau^{(k)}}\,|\, E_0=0]}{1-\beta}
  \label{eq:proof_M_k} \\
  D^{(k)}_\beta(0) &= \EXP \Big[(1-\beta)\sum_{t=0}^{\tau^{(k)}-1} \beta^t d(E_t) + \beta^{\tau^{(k)}} D^{(k)}_\beta(0) \,\Bigm|\, E_0=0 \Big]
  \label{eq:proof_D_k} \\
  N^{(k)}_\beta(0) &= \EXP \Big[ \beta^{\tau^{(k)}}\big( (1-\beta) +
    N^{(k)}_\beta(0)\big) \,\Bigm|\, E_0=0 \Big].
  \label{eq:proof_N_k}
\end{align}
Substituting~\eqref{eq:proof_L_k} and~\eqref{eq:proof_M_k}
in~\eqref{eq:proof_D_k} we get
\begin{equation*}
  D^{(k)}_\beta(0) = (1-\beta) L^{(k)}_\beta(0) 
  + [ 1 - (1 - \beta) M^{(k)}_\beta(0) ] D^{(k)}_\beta(0).
\end{equation*}
Rearranging, we get that
\begin{equation*}
  D^{(k)}_\beta(0) = \frac { L^{(k)}_\beta(0)}{M^{(k)}_\beta(0)}.
\end{equation*}
Similarly, substituting~\eqref{eq:proof_L_k} and~\eqref{eq:proof_M_k} in~\eqref{eq:proof_N_k} we get
\begin{equation*}
  N^{(k)}_\beta(0) = 
  [ 1 - (1 - \beta) M^{(k)}_\beta(0) ] [ (1 - \beta) + N^{(k)}_\beta(0)].
\end{equation*}
Rearranging, we get that
\begin{equation*}
  N^{(k)}_\beta(0) = \frac {1}{M^{(k)}_\beta(0)} - (1-\beta). 
\end{equation*}

The expression for $C^{(k)}_\beta(0;\lambda)$ follows from the definition.

\section{Proofs of results for Model~A}\label{sec:proofs_A}

\subsection{Proof of Theorem~\ref{thm:DIS_AVG}}\label{subsec:proof_DIS_AVG}

\begin{figure}
 \centering
 \includegraphics[width=0.8\linewidth]{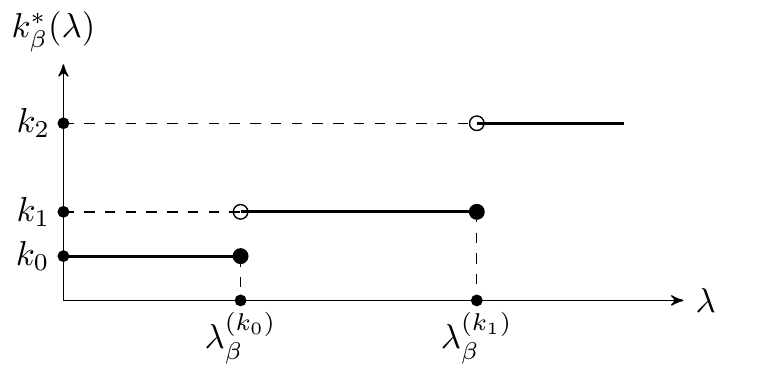}
 \caption{Plot of $k^*_\beta(\lambda)$ for Model~A.}
 \label{fig:k_star_vs_lambda}
\end{figure}

By Proposition~\ref{prop:prelim_c_star_k_star}, $k^*_\beta(\lambda) =
\arg\inf_{k \ge 0} C^{(k)}_\beta(0;\lambda)$ is increasing in $\lambda$. Let
$\mathds K$ denote the set of all possible values of $k^*_\beta(\lambda)$. 
Since $k$ is integer-valued, the plot of $k^*_\beta$ vs $\lambda$ must be a
staircase function as shown in Fig.~\ref{fig:k_star_vs_lambda}. In
particular, there exists an increasing sequence
$\{\lambda^{(k_n)}_\beta\}_{k_n \in \mathds K}$ such that for $\lambda \in
(\lambda^{(k_{n-1})}_\beta, \lambda^{(k_n)}_\beta]$, $k^*_\beta(\lambda) = k_n$. We
will show that for any $k_n$,
\begin{equation}\label{eq:C_at_lambda_k}
  C^{(k_n)}_\beta(0;\lambda^{(k_n)}_\beta) 
  = C^{(k_{n+1})}_\beta(0;\lambda^{(k_n)}_\beta).
\end{equation}
Simplifying \eqref{eq:C_at_lambda_k}, we get that $\lambda^{(k_n)}_\beta$ is
given by~\eqref{eq:lambda}.

\subsubsection*{Proof of \eqref{eq:C_at_lambda_k}}

For any $\lambda \in (\lambda^{(k_{n-1})}_\beta, \lambda^{(k_n)}_\beta]$,
$C^{(k_n)}_\beta(0;\lambda) \le C^{(k_{n+1})}_\beta(0;\lambda)$. In particular, for
$\lambda = \lambda^{(k_n)}_\beta$,
\begin{equation}\label{eq:C_lambdak}
  C^{(k_n)}_\beta(0;\lambda^{(k_n)}_\beta) \le C^{(k_{n+1})}_\beta(0;\lambda^{(k_n)}_\beta).
\end{equation}
Similarly, for any  $\lambda \in (\lambda^{(k_n)}_\beta,
\lambda^{(k_{n+1})}_\beta]$, $C^{(k_{n+1})}_\beta(0;\lambda) \le
C^{(k_n)}_\beta(0;\lambda)$. Since both terms are continuous in $\lambda$, taking
limit as $\lambda \downarrow \lambda^{(k_n)}_\beta$, we get
\begin{equation}\label{eq:C_lambdak_1}
  C^{(k_{n+1})}_\beta(0;\lambda^{(k_n)}_\beta) \le C^{(k_n)}_\beta(0;\lambda^{(k_n)}_\beta).
\end{equation}
Eq.~\eqref{eq:C_at_lambda_k} follows from combining~\eqref{eq:C_lambdak}
and~\eqref{eq:C_lambdak_1}.

\subsubsection{Proof of Part~1)}

By definition of $\lambda^{(k_n)}_\beta$, the strategy $f^{(k_n)}$ is optimal for $\lambda \in (\lambda^{(k_{n-1})}_\beta, \lambda^{(k_n)}_\beta]$.

\subsubsection{Proof of Part~2)}

Recall $C^*_\beta(\lambda) = \inf_{k \ge 0} C^{(k)}_\beta(0; \lambda)$. By
definition, for $\lambda \ge 0$, $C^{(k)}_\beta(0; \lambda)$, is increasing and
affine in $\lambda$. Therefore, its pointwise minimum (over $k$) is
increasing and concave in $\lambda$. 

As shown in part 1), for $\lambda \in (\lambda^{(k_n)}_\beta,
\lambda^{(k_{n+1})}_\beta]$, $C^*_\beta(\lambda) = C^{(k_{n+1})}_\beta(0; \lambda)$,
which is linear (and continuous) in $\lambda$; hence, $C^*_\beta(\lambda)$ is
piecewise linear. Finally, by~\eqref{eq:C_at_lambda_k},
$C^{(k_n)}_\beta(0;\lambda^{(k_n)}) = C^{(k_{n+1})}_\beta(0; \lambda^{(k_n)})$.
Therefore, at the corner points, $\lim_{\lambda \uparrow \lambda^{(k_{n+1})}_\beta}
C^*_\beta(\lambda) = \lim_{\lambda \downarrow \lambda^{(k_{n+1})}_\beta} C^*_\beta(\lambda)$.
Hence, $C^*_\beta(\lambda)$ is continuous in $\lambda$.

\subsection{Proof of Theorem~\ref{thm:constrained_f_star_D}}
\label{subsec:proof_constrained_A}

Note that by definition, $\theta^* \in [0,1]$ and 
\begin{equation}\label{eq:D_theta_star}
  \theta^* N_{\beta}(f^{(k^*)}, g^*) + (1-\theta^*) N_{\beta}(f^{(k^*+1)}, g^*) = \alpha.
\end{equation}

\subsubsection{Proof of Part~1)}
The optimality of $(f^*, g^*)$ relies on the following characterization of
the optimal strategy stated in \cite[Proposition 1.2]{Sennott_avg}. The
characterization was stated for the long-term average setup but a similar
result can be shown for the discounted case as well, for example, by using
the approach of \cite{Borkar}. Also, see \cite[Theorem
8.4.1]{Luenberger1968book} for a similar sufficient condition for general
constrained optimization problem.

A (possibly randomized) strategy $(f^\circ, g^\circ)$ is optimal for a
constrained optimization problem with $\beta \in (0,1]$ if the following
conditions hold:
\begin{enumerate}
  \item[(C1)] $N_\beta(f^\circ,g^\circ) = \alpha$,
  \item[(C2)] There exists a $\lambda^\circ \geq 0$ such that $(f^\circ,g^\circ)$ is optimal for $C_\beta(f,g;\lambda^\circ)$.
\end{enumerate}
We will show that the strategies $(f^*,g^*)$ satisfy (C1) and (C2) with
$\lambda^\circ = \lambda^{(k^*)}_\beta$. 

$(f^*,g^*)$ satisfy (C1) due to~\eqref{eq:D_theta_star}. For $\lambda =
\lambda^{(k^*)}_\beta$, both $f^{(k^*)}$ and $f^{(k^*+1)}$ are optimal for
$C_\beta(f, g; \lambda)$. Hence, any strategy randomizing between them, in
particular $f^*$, is also optimal for $C_\beta(f, g; \lambda)$. Hence $(f^*,
g^*)$ satisfies (C2). Therefore, by \cite[Proposition 1.2]{Sennott_avg},
$(f^*, g^*)$ is optimal for Problem~\ref{prob:P2}.

\subsubsection{Proof of Part~2)}
 
The expression of $k^*$ and $\theta^*$ follow directly from \eqref{eq:k_star}
and \eqref{eq:theta_star}. The form of $D^*_\beta(\alpha)$ given in
(\ref{eq:D_opt_rand}) follows immediately from the fact that $(f^*, g^*)$ is
a Bernoulli randomized simple strategy.

$D^*_\beta(\alpha)$ is the solution to a constrained optimization problem
with the constraint set $\{ (f,g) \, : \, N_\beta(f,g) \le \alpha \}$.
Therefore, it is decreasing and convex in the constraint $\alpha$. The
optimality of $(f^*, g^*)$ implies \eqref{eq:D_opt_rand}. Piecewise linearity
of $D^*_\beta(\alpha)$ follows from \eqref{eq:D_opt_rand}. Finally, by
definition of $\alpha^{(k)}$ and $\theta$, $\lim_{a \uparrow \alpha^{(k)}}
D^*_\beta(\alpha) = D^{(k)}_\beta(0) = \lim_{a \downarrow \alpha^{(k)}}
D^*_\beta(\alpha)$. Hence, $D^*_\beta(\alpha)$ is continuous in $\alpha$.

\section{Proofs of results for Model~B}\label{sec:proofs_B}
\begin{lemma}\label{lem:diff_DNC}
  In Model~B, for $\beta \in (0,1]$,
  \begin{enumerate}
    \item $D^{(k)}_\beta$ and $N^{(k)}_\beta$ are continuous in $k$,
    \item $N^{(k)}_\beta$ is \emph{strictly} decreasing in $k$,
    \item $D^{(k)}_\beta$, $N^{(k)}_\beta$ and $C^{(k)}_\beta$ are
      differentiable in $k$.
  \end{enumerate}
\end{lemma}
\begin{IEEEproof}
  The proof follows from Lemma~\ref{lem:propLM} and
  Theorem~\ref{thm:renewal_DNC}.
\end{IEEEproof}

\subsection{Proof of Theorem~\ref{thm:lambda_k}}

\subsubsection{Proof of Part~1)}

The choice of $\lambda$ implies that $\partial_k C^{(k)}_\beta(0;\lambda)=0$.
Hence strategy $(f^{(k)},g^*)$ is optimal for the given $\lambda$. 

Note that, \eqref{eq:lambda_k} can also be written as $\lambda =
\big((M^{(k)}_\beta(0))^2\partial_k D^{(k)}_\beta(0)\big)/\partial_k
M^{(k)}_\beta(0)$. By Lemma~\ref{lem:propLM}, $\partial_k {M}^{(k)}_\beta(0)
> 0$ and  by Lemma~\ref{lem:DN_lim}, $\partial_k D^{(k)}_\beta(0) \ge 0$.
Hence, for any $k>0$, $\lambda$ given by \eqref{eq:lambda_k} is positive.
This completes the first part of the proof.

\subsubsection{Proof of Part~2)}

The monotonicity and concavity of $C^*_\beta(\lambda)$ follows from the same
argument as in Model~A. 

Note that $k^*_\beta(\lambda) = \arg\inf_{k \ge 0} C^{(k)}_\beta(0;\lambda)$
can take a value $\infty$ (which corresponds to the strategy `never
communicate'). Thus, the domain of $k$ is $\mathds X_{\ge 0}\cup \{\infty\}$,
which is a compact set. Now, $C^*_\beta(\lambda) = \min_{k \in [0,\infty]}
C^{(k)}_\beta(0;\lambda)$, where $C^{(k)}_\beta(0;\lambda)$ is continuous in
both $\lambda$ and $k$. Since, $C^*_\beta(\lambda)$ is pointwise minimum of
bounded continuous functions, where the minimization is over a compact set,
it is continuous.

\subsection{Proof of Theorem~\ref{thm:constrained_cond}}

\subsubsection{Proof of Part~1)}
Recall conditions (C1), (C2), given in
Section~\ref{subsec:proof_constrained_A}, for a strategy to be optimal for a
constrained optimization problem. We will show that for a given $\alpha$,
there exists a $k^*_\beta(\alpha) \in \reals_{\ge 0}$ such that
$(f^{(k^*_\beta(\alpha))}, g^*)$ satisfy conditions (C1) and (C2).

By Lemma~\ref{lem:diff_DNC}, $N^{(k)}_\beta(0)$ is continuous and strictly
decreasing in $k$. It is easy to see that $\lim_{k \to 0} N^{(k)}_\beta(0) =
1$ and $\lim_{k \to \infty} N^{(k)}_\beta(0) = 0$. Hence, for a given $\alpha
\in (0,1)$, there exists a $k^*_\beta(\alpha)$  such that
$N_\beta^{(k^*_\beta(\alpha))}(0) = N_\beta(f^{(k^*_\beta(\alpha))},g^*) =
\alpha$. Thus,  $(f^{(k^*_\beta(\alpha))},g^*)$ satisfies (C1).

Now, for $k^*_\beta(\alpha)$, we can find a $\lambda$ satisfying
\eqref{eq:lambda_k} and hence we have by Theorem~\ref{thm:lambda_k} that
strategy $(f^{(k^*_\beta(\alpha))},g^*)$ is optimal for
$C_\beta(f,g;\lambda)$, and therefore satisfies (C2); and is consequently
optimal for Problem~\ref{prob:P2}.

\subsubsection{Proof of Part~2)}

By Lemma~\ref{lem:diff_DNC}, $\tilde N(k) \DEFINED N^{(k)}_\beta(0)$ is
strictly decreasing and continuous in $k$. Therefore, $\tilde N^{-1}$ exists
and is continuous. Now,
\begin{equation*}
  D^*_\beta(\alpha) = \min_{\{k\,:\, k \le \tilde N^{-1}(\alpha) \}} D^{(k)}_\beta(0),
\end{equation*}
where, by Lemma~\ref{lem:diff_DNC}, $D^{(k)}_\beta(0)$ is continuous in
$k$. Thus, by Berge's maximum theorem, $D^*_\beta(\alpha)$ is continuous in
$\alpha$.

\subsection{Proof of Theorem~\ref{thm:GM_scaling}}
To prove the theorem, we first need to prove the following lemma.

\begin{lemma}\label{lem:scaling_LMDN}
  For Gauss-Markov model (a special case of Model~B), let $L^{(k)}_{\sigma}$
  and $M^{(k)}_{\sigma}$ be the solutions of \eqref{eq:L} and~\eqref{eq:M}
  respectively, when the variance of $W_t$ is $\sigma^2$. Then
  \begin{align}\label{eq:scale_LM}
    L^{(k)}_{\sigma}(e) &= \sigma^2 L^{(k/\sigma)}_{1}\Big(\frac{e}{\sigma}\Big),\quad
    M^{(k)}_{\sigma}(e)=  M^{(k/\sigma)}_{1}\Big(\frac{e}{\sigma}\Big),\\ \label{eq:scale_DN}
    D^{(k)}_{\sigma}(e) &= \sigma^2 D^{(k/\sigma)}_{1}\Big(\frac{e}{\sigma}\Big), \quad
    N^{(k)}_{\sigma}(e) =  N^{(k/\sigma)}_{1}\Big(\frac{e}{\sigma}\Big).
  \end{align}
\end{lemma}

\begin{IEEEproof}
  Define $\hat{L}^{(k)}_{\sigma}(e) \DEFINED \sigma^2 L^{(k/\sigma)}_{1}\Big(\frac{e}{\sigma}\Big)$.
  Now consider,
  \begin{align*}
    [\mathcal B^{(k)}_\sigma \hat L^{(k)}_\sigma](e) &= \int_{-k}^k \phi(n-ae) \hat L^{(k)}_\sigma(n)dn, \quad \forall e \in \reals\\
    &\stackrel{(a)}= \sigma^2 \int_{-k/\sigma}^{k/\sigma} \phi(z-a e/\sigma) L^{(k/\sigma)}_1(z)dz\\
    &= \sigma^2 [\mathcal B^{(k/\sigma)}_1 L^{(k/\sigma)}_1] (e/\sigma),
  \end{align*}
  where $(a)$ uses a change of variables $n=\sigma z$. Therefore,
  \begin{align*}
    \Big[\hat{L}^{(k)}_{\sigma} - \beta \mathcal{B}^{(k)}_{\sigma} \hat{L}^{(k)}_{\sigma}\Big](e) &=
    \sigma^2 \Big[L^{(k/\sigma)}_{1} - \beta \mathcal{B}^{(k/\sigma)}_{1} L^{(k/\sigma)}_{1}\Big]\Big(\frac{e}{\sigma}\Big)\\
    &= \sigma^2 \frac{e^2}{\sigma^2} =  e^2.
  \end{align*}
  But, by Lemma~\ref{lem:propLM}, the above equation has a unique solution
  ${L}^{(k)}_\sigma$. Therefore $L^{(k)}_\sigma = \hat{L}^{(k)}_{\sigma}$.

  A similar argument may be used to prove the scaling of $M^{(k)}_\sigma$.
  The scaling of $D^{(k)}_\sigma$ and $N^{(k)}_\sigma$ follow from
  Theorem~\ref{thm:renewal_DNC}.
\end{IEEEproof}

\subsubsection*{Proof of Theorem~\ref{thm:GM_scaling}}
The theorem follows from Lemma~\ref{lem:scaling_LMDN},
Theorem~\ref{thm:renewal_DNC} and elementary algebra.

\section{Proofs of results for Example~\ref{ex:BD}}\label{sec:proofs_example}
\begin{lemma} \label{lem:top}
  Define for $\beta \in (0,1]$
  \begin{gather*}
  K_\beta = -2 - \frac{(1-\beta)}{\beta p} \quad \text{and} \quad m_\beta = \cosh^{-1}(-K_\beta/2)
  \end{gather*}
 Then,
  \begin{equation*}
    [Q^{(k)}_\beta]_{ij} = \frac 1 {\beta p} \frac { [A^{(k)}_\beta]_{ij} } {
      b^{(k)}_\beta },
      \quad i,j \in S^{(k)},
  \end{equation*}
  where, for $\beta \in (0,1)$,
  \begin{align*}
    [A^{(k)}_\beta]_{ij} &= \cosh( (2k - |i-j|)m_\beta ) - \cosh( (i+j) m_\beta), \\
    b^{(k)}_\beta &= \sinh(m_\beta) \sinh(2km_\beta);
  \end{align*}
  and for $\beta = 1$,
  \begin{align*}
    [A^{(k)}_1]_{ij} &= ( k - \max\{i, j\})(k + \min\{i,j\}), \\
    b^{(k)}_1 &= 2k. 
  \end{align*}
  In particular, the elements $[Q^{(k)}_\beta]_{0j}$ are given as follows. For
  $\beta \in (0,1)$,
  \begin{equation}\label{Q_beta_BD}
    [Q^{(k)}_\beta]_{0j} = \frac 1 {\beta p}
    \frac {\cosh( (2k - |j|) m_\beta) - \cosh(j m_\beta)}
          {2 \sinh(m_\beta) \sinh( 2k m_\beta) },
  \end{equation}
  and for $\beta = 1$,
  \begin{equation}\label{Q_1_BD}
  [Q^{(k)}_1]_{0j} = \frac{k-|j| }{ 2 p }.
  \end{equation}
\end{lemma}
\begin{IEEEproof}
  The matrix $I_{2k-1} - \beta P^{(k)}$ is a symmetric tridiagonal matrix
  given by
  \begin{equation*}
    I_{2k-1} - \beta P^{(k)} = - \beta p 
    \begin{bmatrix}
      K_\beta & 1       & 0       & \cdots & \cdots & 0 \\
      1       & K_\beta & 1       & 0      & \cdots & 0 \\
      0       & 1       & K_\beta & 1      & \cdots & 0 \\
      \vdots  & \ddots  & \ddots  & \ddots & \ddots & \vdots \\
      0       & \cdots  & 0       & 1      & K_\beta& 1 \\
      0       & 0       & \cdots  & 0       & 1     & K_\beta 
    \end{bmatrix}.
  \end{equation*}
  $Q^{(k)}_\beta$ is the inverse of the above matrix. The inverse of the
  tridiagonal matrix in the above form with $K_\beta \le -2$ are computed in
  closed form in~\cite{HuOConnell:1996}. The result of the lemma follows from
  these results.   
\end{IEEEproof}

\subsection{Proof of Lemma~\ref{lem:DN_BDMC}}
By substituting the expression for $Q^{(k)}_\beta$ from Lemma~\ref{lem:top}
in the expressions for $L^{(k)}_\beta$ and $M^{(k)}_\beta$ from
Proposition~\ref{prop:expressLM_A}, we get that
\begin{enumerate}
  \item For $\beta \in (0,1)$, 
    \begin{align*}
      L^{(k)}_\beta(0) &= \frac { \sinh(k m_\beta) - k \sinh (m_\beta) }
      { 4 \beta p \sinh^2(m_\beta/2) \sinh(m_\beta) \cosh(k m_\beta)} ,
      \\
      M^{(k)}_\beta(0) &= \frac { \sinh^2 (k m_\beta /2 ) }
      { 2 \beta p \sinh^2(m_\beta/2) \cosh(k m_\beta)}.
    \end{align*}

  \item For $\beta = 1$,
    \begin{equation*}
      L^{(k)}_1(0) = k(k^2 - 1)/ (6p), \quad
      M^{(k)}_1(0) = k^2 / (2p).
    \end{equation*}
\end{enumerate}
The results of the lemma follow using the above expressions and Theorem~\ref{thm:renewal_DNC}. The expression for $\lambda^{(k)}_1$ is obtained by plugging the expressions of $D^{(k+1)}_1$, $D^{(k)}_1$, $N^{(k+1)}_1$, and $N^{(k)}_1$ in~\eqref{eq:lambda}.

\section{Conclusion}\label{sec:conclusion}
We characterize two fundamental limits of remote estimation of autoregressive
Markov processes under communication constraints. First, when each
transmission is costly, we characterize the minimum achievable cost of
communication plus estimation error. Second, when there is a constraint on
the average number of transmissions, we characterize the minimum achievable
estimation error. 

We also identify transmission and estimation strategies that achieve these
fundamental limits. The structure of these optimal strategies had been
previously identified by using dynamic programming for decentralized
stochastic control systems. In particular, the optimal transmission strategy
is to transmit when the estimation error process exceeds a threshold and the
optimal estimation strategy is to select the transmitted state as the
estimate, whenever there is a transmission. We use ideas based on renewal
theory to identify the performance of a generic strategy that has such a
structure. For the case of costly communication, we identify the value of
communication cost for which a particular threshold-based strategy is
optimal; for the case of constrained communication, we identify (possibly
randomized) threshold-based strategies that achieve the communication
constraint. 

These results are derived under idealized assumptions on the communication
channel: communication is noiseless and without any constraint on the
transmission rate or the transmission bandwidth. Under these assumptions, the
error process resets after each transmission (see Remark~\ref{rem:reset}).
This reset property is critical to derive the structure of optimal
transmission and estimation strategies (Theorems~\ref{thm:inf_hor}
and~\ref{thm:fin_hor}). In the absence of such a structural result, the
solution methodology developed in this paper does not work and the optimal
transmission and estimation strategies have to be identified by numerically
solving the (decentralized) dynamic programs described in
\cite{WalrandVaraiya:1983,MT:real-time}. 

Having said that, the transmission and estimation strategies described in
Theorems~\ref{thm:inf_hor} and~\ref{thm:fin_hor} may be used as heuristic
sub-optimal strategies when the communication channel does not satisfy the
idealized assumptions described above. In that case, it may be possible to
use the solution methodology developed in this paper to obtain performance
bounds on such strategies. 

A similar remark holds for multi-dimensional autoregressive processes. It is
reasonable to expect (although we are not aware of a proof of this statement)
that for multi-dimensional autoregressive processes, the optimal estimation
strategy will be similar to that described in Theorems~\ref{thm:inf_hor}
and~\ref{thm:fin_hor} while the optimal transmission strategy will be to
transmit when the error process lies outside a (multi-dimensional) ellipsoid.
The performance of such strategies can be evaluated using the solution
methodology developed in this paper. The renewal relationships derived in
Theorem~\ref{thm:renewal_DNC} also hold for multi-dimensional autoregressive
processes. The only difference is that $L^{(k)}_\beta(0)$ and
$M^{(k)}_\beta(0)$ are computed by solving multi-dimensional Fredholm
integral equations of the second kind. The optimal transmission strategies
can then be computed by solving multi-dimensional versions of
\eqref{eq:lambda_k} (for costly communication) and \eqref{eq:N_star} (for
constrained communication). However, it is not immediately clear whether
these equations will have a unique solution. Further investigation is
required to obtain algorithms that identify the optimal transmission
ellipsoid.

Finally, the solution methodology developed in this paper to identify optimal
thresholds is also of independent interest. In various applications of Markov
decision processes threshold strategies are optimal. The approach developed
in this paper is directly applicable to such models.
  
\section*{Acknowledgments}

The authors are grateful to M.~Madiman, A.~Molin, A.~Paranjape,
V.~Subramanian, and S.~Y\"uksel for useful discussions.

\appendices

\section{Proof of Lemma~\ref{lem:propLM}}\label{lem:proof_propLM}
Let $\|\cdot\|_\infty$ denote the sup-norm, i.e., for any $v: S^{(k)} \to
\reals$,
\begin{equation*}
  \|v\|_\infty = \sup_{e \in S^{(k)}} |v(e)|.
\end{equation*}
To prove the lemma, let us first prove the following:
\begin{lemma}\label{lem:contraction}
  For $\beta \in (0,1)$, for both Models~A and~B, the operator $\beta
  \mathcal{B}^{(k)}$ is a contraction, i.e., for any $v: S^{(k)} \to \reals$,
  \begin{equation*}
    \| \beta \mathcal{B}^{(k)}v\|_\infty \le \beta \|v\|_\infty.
  \end{equation*}
  Thus, for any bounded $h: S^{(k)} \to \reals$, the equation 
  \begin{equation}\label{eq:op_eqn}
    v = h+\beta \mathcal{B}^{(k)}v
  \end{equation}
  has a unique bounded solution $v$. In addition, if $h$ is continuous, then
  $v$ is continuous.
\end{lemma}

\begin{IEEEproof}
  We state the proof for Model B. The proof for Model A is similar. By the
  definition of sup-norm, we have that for any bounded $v$
  \begin{align*}
    \|\beta \mathcal{B}^{(k)} v\|_\infty &= \beta \sup_{e \in (-k,k)}
    \int_{-k}^k \phi(w-ae) v(w) dw\\
    & \le \beta \sup_{e \in (-k,k)} \|v\|_\infty \int_{-k}^k \phi(w-ae) dw\\
    & \le \beta \|v\|_\infty, \quad (\textrm{since $\phi$ is a pdf}).
  \end{align*}
  Hence, $\beta \mathcal{B}^{(k)}$ is a contraction.

  Now, consider the operator $\mathcal{B}'$ given as: $\mathcal{B}' v = h +
  \beta \mathcal{B}^{(k)} v$. Then we have,
  \begin{equation*}
    \|\mathcal{B}' (v_1-v_2)\|_\infty 
    = \beta  \|\mathcal{B}^{(k)} (v_1-v_2)\|_\infty
    \le \beta \|v_1-v_2\|_\infty.
  \end{equation*}
  Since $\beta \in (0,1)$ and the space of bounded real-valued functions is
  complete, by Banach fixed point theorem, $\mathcal{B}'$ has a unique fixed
  point.

  If $h$ is continuous, we can define $\mathcal{B}^{(k)}$ and $\mathcal{B}'$
  as operators on the space of continuous and bounded real-valued function
  (which is complete). Hence, the continuity of the fixed point follows also
  from Banach fixed point theorem.
\end{IEEEproof}

\subsection*{Proof of Lemma~\ref{lem:propLM}}

The solutions of equations~\eqref{eq:L} and~\eqref{eq:M} exist due to
Lemma~\ref{lem:contraction}. 
\begin{enumerate}
  \item[(a)] Consider $k$, $l \in \mathds X_{\ge 0}$ such that $k<l$. A sample path starting from $e \in S^{(k)}$ must escape $S^{(k)}$ before it escapes $S^{(l)}$. Thus $L^{(l)}_\beta(e) \ge L^{(k)}_\beta(e)$. In addition, the above inequality is strict because $W_t$ has a unimodal distribution. Similar argument holds for $M^{(k)}_\beta$.
  \item[(b)] The continuity and differentiability can be proved from elementary algebra. 
    See the supplementary material for details.
  \item[(c)] The limit holds since $L^{(k)}_\beta(e)$ and $M^{(k)}_\beta(e)$ are continuous functions of $\beta$.
\end{enumerate}

\section{Proof of Proposition~\ref{prop:prelim_c_star_k_star}}
\label{app:prelim_c_star_k_star}
\begin{enumerate}
  \item $C^{(l)}_\beta(0;\lambda) - C^{(k)}_\beta(0;\lambda)=
    (D^{(l)}_\beta(0) - D^{(k)}_\beta(0)) - \lambda (N^{(k)}_\beta(0) -
    N^{(l)}_\beta(0))$. By Lemma~\ref{lem:propLM} and
    Theorem~\ref{thm:renewal_DNC}, $N^{(k)}_\beta(0) - N^{(l)}_\beta(0)$ is
    positive, hence $C^{(l)}_\beta(0;\lambda) - C^{(k)}_\beta(0;\lambda)$ is
    decreasing in $\lambda$. Hence $C^{(k)}_\beta(0;\lambda)$ is submodular.
  \item Note that $k^*_\beta(\lambda) = \arg\inf_{k \ge 0}
    C^{(k)}_\beta(0;\lambda)$ can take a value $\infty$ (which corresponds to
    the strategy `never communicate'). Thus, the domain of $k$ is $\mathds
    X_{\ge 0}\cup \{\infty\}$, which is compact. Hence, by~\cite[Theorem
    2.8.2]{Topkis}, $k^*_\beta$ is increasing in $\lambda$.
\end{enumerate}

\section{Proofs of Propositions~\ref{prop:V_plus_minus} and~\ref{prop:EI}}
\label{app:V_plus_minus_EI}

We prove the results for Model~A when the horizon $T$ is finite. The results
then follow by taking limits as $T \rightarrow \infty$. The proofs for Model~B
are almost identical.

The value function for the finite horizon setup for $\beta \in (0,1]$is given by $V_{\beta,T+1} =
0$ and for $t=T,\cdots,1$
\begin{align}
  V_{\beta,t} (e;\lambda) &= \min \Big \{(1-\beta)\lambda + \beta
    \sum_{n=-\infty}^\infty p_n V_{\beta,t+1}(n;\lambda),\notag \\ 
  & \qquad  (1-\beta) d(e) + \beta \sum_{n=-\infty}^\infty p_{n-ae} V_{\beta,t+1}(n;\lambda)\Big\}.
\label{eq:DP_fin_hor}
\end{align}
The value functions $V^{(+)}_t$ and $V^{(-)}_t$ are defined similarly. 

For ease of notation, we drop $\beta$ and $\lambda$ in the rest of the
discussion in this Appendix.

\begin{lemma} \label{lem:even}
  The value functions $V_t(\cdot)$, $V^{(+)}_t(\cdot)$ and $V^{(-)}_t(\cdot)$ are even.
\end{lemma}

\begin{IEEEproof}
  For all $a \in \mathds X$, the per-step costs $d(e)$ and $\lambda$ are even
  and the transition probabilities $P_{en}(0) = p_{n - ae}$ and $P_{en}(1) = p_n$
  satisfy $P_{en}(u)  = P_{(-e)(-n)}(u)$ for $u \in \{0,1\}$. 
  Therefore, $V_t(e)$ is even~\cite[Theorem~1]{CM:even-increasing-icc}. A
  similar argument holds for $V^{(+)}_t(e)$ and $V^{(-)}_t(e)$. 
\end{IEEEproof}

\begin{lemma}\label{lem:V_plus_minus_equal}
  For the finite horizon setup, $V^{(+)}_t(e) = V^{(-)}_t(e)$.
\end{lemma}

\begin{IEEEproof}
  We prove the result by backward induction. The result is trivially true for
  $T+1$ as $V^{(+)}_{T+1}(e) = V^{(-)}_{T+1}(e) = 0$, which forms the basis of
  the induction. Assume $V^{(+)}_{t+1}(e) = V^{(-)}_{t+1}(e)$ for all $e \in
  \mathds X$.
  Define
  \begin{equation*}
    \hat V^{(+)}_t(e) = \SUMN p_{n-ae} V^{(+)}_{t+1}(n), \quad \hat V^{(-)}_t(e) = \SUMN p_{n+ae} V^{(-)}_{t+1}(n).
  \end{equation*}
  Then
  \begin{align*}
    \hat V^{(+)}_t(e) &= \sum_{n = -\infty}^{\infty} p_{n-ae} V^{(+)}_{t+1}(n) = \sum_{-n = -\infty}^{\infty} p_{-n-ae} V^{(+)}_{t+1}(-n) \\
    &\stackrel{(a)}= \SUMN p_{n+ae} V^{(+)}_{t+1}(n)
    \stackrel{(b)}= \SUMN p_{n+ae} V^{(-)}_{t+1}(n) = \hat V^{(-)}_t(e),
  \end{align*}
  where $(a)$ uses $p$ and $V^{(+)}_{t+1}$ are even and $(b)$ uses the
  induction hypothesis. Substituting this back in the definition of
  $V^{(+)}_t(e)$ and $V^{(-)}_t(e)$, we get that $V^{(+)}_t(e) = V^{(-)}_t(e)$.
  Therefore, the result is true by induction.
\end{IEEEproof}

\begin{lemma} \label{lem:Q}
  For $m,e \in \mathds X_{\ge 0}$, define
  \begin{equation*}
    Q(m|e,0) = \sum_{n : |n| \ge m} p_{n - ae}
    \quad\text{and}\quad
    Q(m|e,1) = \sum_{n : |n| \ge m} p_n.
  \end{equation*}
  Then, for all $e, m \in \mathds X_{\ge 0}$ and $a > 0$,  $Q(m|e,0)$ and $Q(m|e,1)$ are
  increasing in~$e$. 
\end{lemma}

We will prove this Lemma later. 

\begin{definition}
  A function $f \colon \mathds X \to \reals$ is called even and increasing
  on $\mathds X_{\ge 0}$  if for all $x \in \mathds
  X_{\ge 0}$, $f(x) = f(-x)$ and $f(x) \le f(x+1)$. 
\end{definition}

\begin{lemma} \label{lem:EI}
  The value function $V_t(e)$ is even and increasing on $\mathds X_{\ge 0}$. 
\end{lemma}
\begin{IEEEproof}
  We have already shown that $V_t(e)$ is even. For $a > 0$, the properties
  described in the proof of Lemma~\ref{lem:even} and the statement
  Lemma~\ref{lem:Q} imply that $V_t(e)$ is
  even and increasing~\cite[Theorem~1]{CM:even-increasing-icc}. Now,
  Lemma~\ref{lem:V_plus_minus_equal} implies that $V_t(e)$ is also even and
  increasing for $a < 0$. 
\end{IEEEproof}

\begin{IEEEproof}[Proofs of Propositions~\ref{prop:V_plus_minus} and~\ref{prop:EI}]  
  The result follows from Lemmas~\ref{lem:V_plus_minus_equal} and~\ref{lem:EI}
  by taking the limit $T \rightarrow \infty$, since equality is preserved
  under limits.
\end{IEEEproof}

\begin{IEEEproof}[Proof of Lemma~\ref{lem:Q}]
  $Q(m|e,1)$ is independent of~$e$. Define $R(m|e) = \sum_{n : |n| \le m }
  p_{n - e}$. Then, $Q(m|e,0) = 1 - R(m|ae)$. To show $Q(m|e,0)$ is increasing
  in~$e$, it suffices to show that $R(m|ae) \ge R(m|ae+1)$ (which implies that
  $R(m|ae) \ge R(m|ae+a)$).

  Now consider
  \begin{equation*}
    R(m|ae) - R(m|ae+1) = p_{m-ae} - p_{-m-ae-1} = p_{m-ae} - p_{m+ae+1}.
  \end{equation*}
  If $m \ge ae$, then $0 \le m - ae < m + ae + 1$, hence, $p_{m-ae} \ge
  p_{m+ae+1}$. If $m < ae$, then $0 < ae - m < m + ae + 1$, hence $p_{m-ae}
  = p_{ae -m} \ge p_{m+ae+1}$. Thus, in both cases, $R(m|ae) \ge R(m|ae+1)$. 
\end{IEEEproof}

\section{Proof of Part~3) of Lemma~\ref{lem:DN_lim}}\label{app:D_inc}
By Lemma~\ref{lem:propLM}, $M^{(k)}_\beta(e)$ is strictly increasing in~$k$;
therefore, by Theorem~\ref{thm:renewal_DNC}, $N^{(k)}_\beta(e)$ is strictly
decreasing in~$k$.

We prove the monotonicity of $D^{(k)}_\beta$ in $k$ for Model A for
$\beta \in (0,1)$. The result for $\beta = 1$ follows by taking limit $\beta
\uparrow 1$. The result for Model~B is similar. Based on
Lemma~\ref{lem:V_plus_minus_equal}, we restrict attention to $a > 0$.

For any $\beta \in (0,1)$ and $k \in \integers_{\ge 0}$, define the operator
$\mathcal T^{(k)}: (\integers \rightarrow \reals) \rightarrow (\integers \rightarrow
\reals)$ as follows. For any $D: \integers \rightarrow \reals$,
\begin{equation}\label{eq:TD}
  [\mathcal T^{(k)}D](e) = 
  \begin{cases}
    \beta [ \mathcal B D ] (0), & \text{if $|e| \ge k$} \\
    (1-\beta)d(e)+\beta [ \mathcal B D](e)  & \text{if $|e| < k$}.
  \end{cases}
\end{equation}
This operator is the Bellman operator for evaluating strategy $f^{(k)}$.
Hence, it is a contraction and $D^{(k)}$ is the unique fixed point of
$\mathcal T^{(k)}$. 

Define $D^{(k,0)}_\beta =  D^{(k)}_\beta$, and for $m \in \integers_{> 0}$,
$D^{(k,m)}_\beta = \mathcal T^{(k+1)}  D^{(k,m-1)}_\beta$. 

From Lemma~\ref{lem:Q}
and~\cite[Lemma~2]{CM:even-increasing-icc}, we get that for any $e \in
\integers_{\ge 0}$, 
\begin{equation*}
  \sum_{n =-\infty}^\infty p_{n-ae} D^{(k)}_\beta (n) \ge
  \sum_{n =-\infty}^\infty p_n D^{(k)}_\beta (n),
\end{equation*}
or equivalently, $[\mathcal B D^{(k)}_\beta](e) \ge [\mathcal B
D^{(k)}_\beta](0)$.

For $|e| = k$, $D^{(k,1)}_\beta(e) = (1-\beta) d(e) + \beta [\mathcal B
D^{(k)}_\beta](e)$ and $D^{(k)}_\beta(e) = \beta [\mathcal B
D^{(k)}_\beta](0)$; hence, $D^{(k,1)}_\beta(e) > D^{(k)}_\beta(e)$. For
$|e| \neq k$, $D^{(k,1)}_\beta(e) = D^{(k)}_\beta(e)$ because both terms
have the same expression. Hence, for all $e \in \integers$, 
\begin{equation*}
  D^{(k,1)}_\beta(e) \ge D^{(k)}_\beta (e), 
  \quad\text{or}\quad
  D^{(k,1)}_\beta \ge D^{(k)}_\beta.
\end{equation*}
If we apply the operator $\mathcal T^{(k+1)}$ to both sides, the monotonicity
of $\mathcal T^{(k+1)}$ implies that $D^{(k,2)}_\beta \ge D^{(k,1)}_\beta \ge
D^{(k)}_\beta$. Proceeding this way, we get that for any $m > 0$,
\begin{equation}\label{eq:Dm}
  D^{(k+m)}_\beta \ge D^{(k)}_\beta.
\end{equation}
Note that $\lim_{m \to \infty} D^{(k+m)}_\beta = D^{(k+1)}_\beta$, because
$D^{(k+1)}_\beta$ is the unique fixed point of the operator $\mathcal
T^{(k+1)}$. Thus, taking limit $m \to \infty$ in~\eqref{eq:Dm}, we get that
$D^{(k+1)}_\beta \ge D^{(k)}_\beta$.

% Generated by IEEEtran.bst, version: 1.14 (2015/08/26)

\begin{IEEEbiography}[{\includegraphics[width=1in]{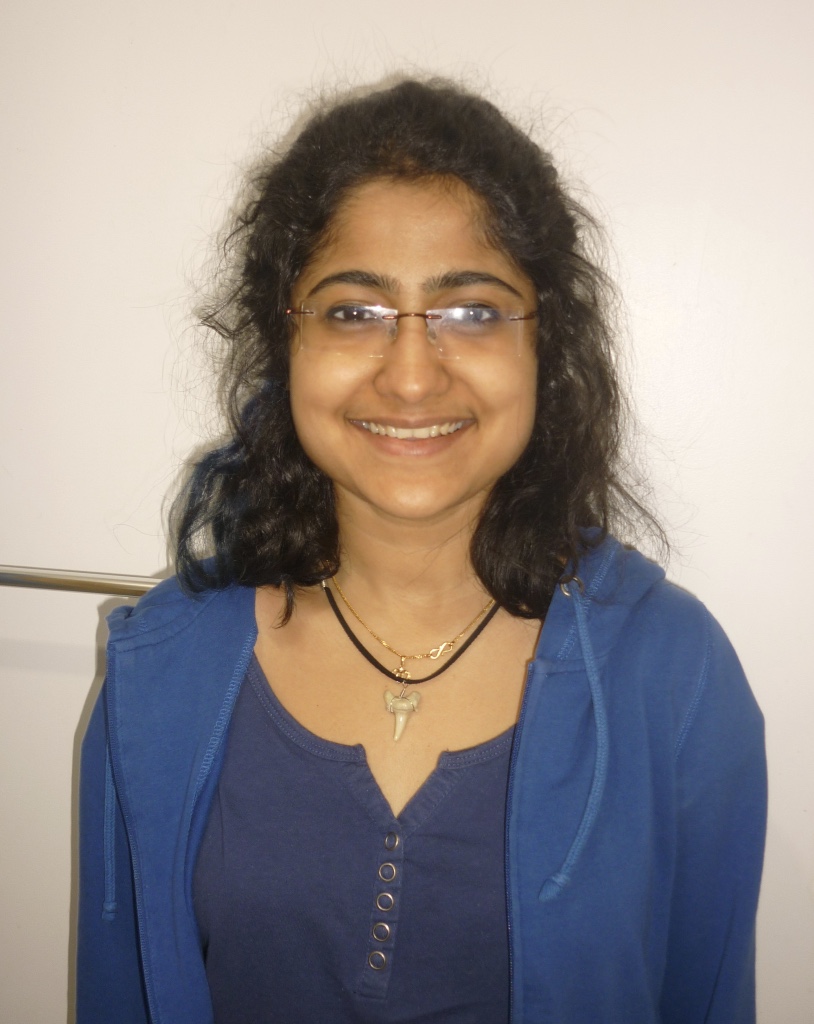}}]
{Jhelum Chakravorty} is a doctoral student of Electrical and Computer Engineering at McGill University, Montreal, Canada. She received the B.E. degree in Electrical Engineering from Jadavpur University, India in 2007 and M.Tech degree in Systems and Control Engineering from Indian Institute of Technology, Bombay in 2010. Prior to the doctoral study, she worked as a research assistant at Indian Institute of Technology, Bombay during 2007--2008, and as a research associate at Indian Institute of Science during 2010--2011. Her current area of research is decentralized stochastic control, information theory and real-time communication.  
\end{IEEEbiography}

\begin{IEEEbiography}[{\includegraphics[width=1in]{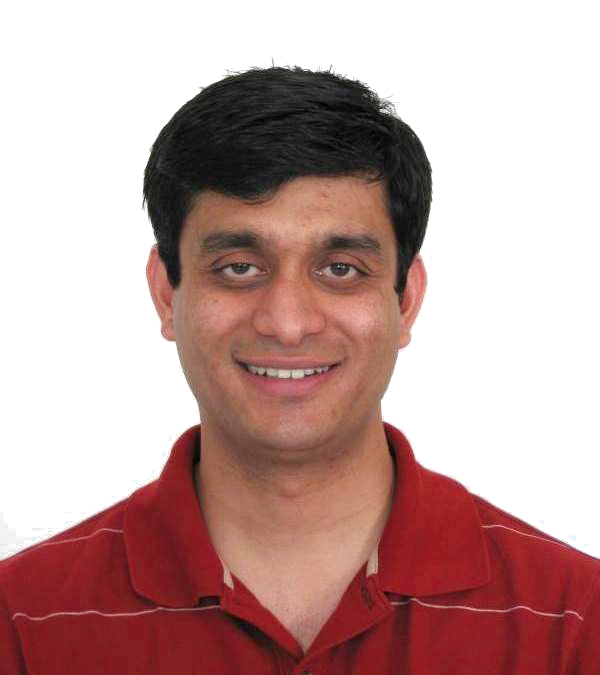}}]
  {Aditya Mahajan} (S'06--M'09--SM-14) is Associate Professor of Electrical
  and Computer Engineering at McGill University, Montreal, Canada. He is 
  member of the McGill Center of Intelligent Machines (CIM),  member of
  Groupe d'\'etudes et de recherche en analyse des d\'ecisions (GERAD), 
  senior member of the IEEE, member of SIAM, and member of Professional
  Engineers Ontario (PEO).
  He currently serves as an Associate Editor of the IEEE Control Systems
  Society Conference Editorial Board. 
  
  He received B.Tech degree in Electrical Engineering from the Indian
  Institute of Technology, Kanpur, India in 2003 and MS and PhD degrees
  in Electrical Engineering and Computer Science from the University of
  Michigan, Ann Arbor, USA in 2006 and 2008. From 2008 to 2010, he was
  postdoctoral researcher in the department of Electrical Engineering at Yale
  University, New Haven, CT, USA. Since 2010, he is with the department of
  Electrical and Computer Engineering at McGill University, Montreal, Canada.
  His principal research interests include decentralized stochastic control,
  team theory, multi-armed bandits, real-time communication, information
  theory, and discrete event systems.
\end{IEEEbiography}

\clearpage

\title {Supplementary material for ``Fundamental limits of remote estimation
of autoregressive Markov processes under communication constraints''}

\author{Jhelum Chakravorty and Aditya Mahajan}

\maketitle

\section{Proof of the structural results} \label{app:structure}

The results of~\cite{NayyarBasarTeneketzisVeeravalli:2013} relied on the notion
of ASU (almost symmetric and unimodal) distributions introduced
in~\cite{HajekMitzelYang:2008}. 
\begin{definition}[Almost symmetric and unimodal distribution]
 A probability distribution $\mu$ on $\integers$ is almost symmetric and
 unimodal (ASU) about a point $a \in \integers$ if for every $n \in
 \integers_{\ge 0}$,
 \[
   \mu_{a + n} \ge \mu_{a - n} \ge \mu_{a + n + 1}.
 \]
\end{definition}

A probability distribution that is ASU around~$0$ and even (i.e., $\mu_n =
\mu_{-n}$) is called ASU and even. Note that the definition of ASU and even is
equivalent to even and decreasing on~$\integers_{\ge 0}$. 

\begin{definition}[ASU Rearrangement]
 The ASU rearrangement of a probability distribution $\mu$, denoted by
 $\mu^{+}$, is a permutation of $\mu$ such that for every $n \in
 \integers_{\ge 0}$,
 \[
   \mu^{+}_n \ge \mu^{+}_{-n} \ge \mu^{+}_{n+1}.
 \]
\end{definition}

We now introduce the notion of majorization for distributions supported
over~$\integers$, as defined in~\cite{WangWooMadiman:2014}.

\begin{definition}[Majorization]
 Let $\mu$ and $\nu$ be two probability distributions defined over
 $\integers$. Then $\mu$ is said to majorize $\nu$, which is denoted by $\mu
 \succeq_m \nu$, if for all $n \in \integers_{\ge 0}$,
 \begin{align*}
   \sum_{i = -n}^n \mu^{+}_i &\ge \sum_{i=-n}^n \nu^{+}_i, \\
   \sum_{i = -n}^{n+1} \mu^{+}_i &\ge \sum_{i=-n}^{n+1} \nu^{+}_i. 
 \end{align*}
\end{definition}

The structure of optimal estimator in Theorem~8 were proved in two steps
in~\cite{NayyarBasarTeneketzisVeeravalli:2013}. The first step relied on the
following two results.
\begin{lemma}
 \label{lem:ASU}
 Let $\mu$ and $\nu$ be probability distributions with finite support defined
 over $\integers$. If $\mu$ is ASU and even and $\nu$ is ASU about $a$, then
 the convolution $\mu \ast \nu$ is ASU about $a$.
\end{lemma}
\begin{lemma}
 \label{lem:ASU-M}
 Let $\mu$, $\nu$, and $\xi$ be probability distributions with finite
 support defined over $\integers$. If $\mu$ is ASU and even, $\nu$ is ASU, and
 $\xi$ is arbitrary, then $\nu \succeq_m \xi$ implies that $\mu \ast \nu
 \succeq_m \mu \ast \xi$.
\end{lemma}
These results were originally proved in~\cite{HajekMitzelYang:2008} and were
stated as Lemmas~5 and~6 in~\cite{NayyarBasarTeneketzisVeeravalli:2013}.

The second step (in the proof of structure of optimal estimator in Theorem~8)
in~\cite{NayyarBasarTeneketzisVeeravalli:2013} relied on the following result.
\begin{lemma}
 \label{lem:ineq}
 Let $\mu$ be a probability distribution with finite support defined over
 $\integers$ and $f \colon \integers \to \reals_{\ge 0}$. Then,
 \[
   \SUMN f(n) \mu_n \le \SUMN f^{+}(n) \mu^{+}_n.
 \]
\end{lemma}

We generalize the results of Lemmas~\ref{lem:ASU}, \ref{lem:ASU-M}, and
\ref{lem:ineq} to distributions over $\integers$ with possibly countable support. With these
generalizations, we can follow the same two-step approach
of~\cite{NayyarBasarTeneketzisVeeravalli:2013} to prove the structure of optimal estimator as given in 
Theorem~8.

The structure of optimal transmitter in Theorem~8
in~\cite{NayyarBasarTeneketzisVeeravalli:2013} only relied on the structure of optimal estimator. The exact same proof works in our model as well.

\subsection {Generalization of Lemma~\ref{lem:ASU} to distributions supported over~$\integers$}

The proof argument is similar to that presented in~\cite[Lemma 6.2]{HajekMitzelYang:2008}.
We first prove the results for $a = 0$. Assume that $\nu$ is ASU and even. For any $n \in \integers_{\ge 0}$, let $r^{(n)}$ denote the rectangular function from $-n$ to $n$, i.e.,
\[
r^{(n)}(e) = \begin{cases} 
    1, & \text{if $|e| \le n$}, \\
    0, & \text{otherwise}.
  \end{cases}
\]

Note that any ASU and even distribution $\mu$ may be written as a sum of rectangular functions as follows:
\[
  \mu = \sum_{n=0}^\infty (\mu_n - \mu_{n+1}) r^{(n)}.
\]
It should be noted that $\mu_n - \mu_{n+1} \ge 0$ because $\mu$ is ASU and even. $\nu$ may also be written in a similar form.

The convolution of any two rectangular functions $r^{(n)}$ and $r^{(m)}$ is ASU and even. Therefore, by the distributive property of convolution, the convolution of $\mu$ and $\nu$ is also ASU and even.

The proof for the general $a \in \integers$ follows from the following facts:
\begin{enumerate}
 \item Shifting a distribution is equivalent to convolution with a shifted
   delta function.
 \item Convolution is commutative and associative.
\end{enumerate}

\subsection {Generalization of Lemma~\ref{lem:ASU-M} to distributions supported over~$\integers$}

We follow the proof idea of~\cite[Theorem II.1]{WangWooMadiman:2014}. For any
probability distribution $\mu$, we can find distinct indices $i_j$, $|j| \le n$
such that $\mu(i_j)$, $|j| \le n$, are the $2n + 1$ largest values of $\mu$. Define
\[
 \mu_n(i_j) = \mu(i_j),
\]
for $|j| \le n$ and $0$ otherwise. Clearly, $\mu_n \uparrow \mu$ and if $\mu$
is ASU and even, so is $\mu_n$. 

Now consider the distributions $\mu$, $\nu$, and $\xi$ from Lemma~\ref{lem:ASU-M}
but without the restriction that they have finite support. For every $n \in
\integers_{\ge 0}$, define $\mu_n$, $\nu_n$, and $\xi_n$ as above. Note that
all distributions have finite support and $\mu_n$ is ASU and even and $\nu_n$ is
ASU. Furthermore, since the definition of majorization remain unaffected by
truncation described above, $\nu_n \succeq_m \xi_n$.
Therefore, by Lemma~\ref{lem:ASU-M}, 
\[
 \mu_n \ast \nu_n \succeq_m \mu_n \ast \xi_n.
\]
By taking limit over $n$ and using the monotone convergence theorem, we get
\[
 \mu \ast \nu \succeq_m \mu \ast \xi.
\]

\subsection {Generalization of Lemma~\ref{lem:ineq} to distributions supported over~$\integers$}

This is an immediate consequence of~\cite[Theorem II.1]{WangWooMadiman:2014}.

\section{Proof of (b) of Lemma~1}\label{app:cont_diff_LM}

Note that for any bounded $v$, $\|\mathcal{B}^{(k)}v\|_\infty$ is bounded and increasing in $k$. We show that $L^{(k)}_\beta(e)$ is continuous and differentiable in $k$. Similar argument holds for $M^{(k)}_\beta(e)$.

We show the differentiability in $k$. Continuity follows from the fact that differentiable functions are continuous. 
Note that $L^{(k)}_\beta(e)$ and $M^{(k)}_\beta(e)$ are even functions of $e$. Now, for any $\varepsilon >0$ we have
\begin{align*}
& L^{(k+\varepsilon)}_\beta (e) - L^{(k)}_\beta (e) \\
&= \beta \int_{-k}^k \phi(w-ae) [L^{(k+\varepsilon)}_\beta (w)-L^{(k)}_\beta (w)] dw \\
& \hskip 2em  + 2 \beta \int_k^{k+\varepsilon} \phi(w-ae) L^{(k+\varepsilon)}_\beta (w) dw\\
& = \beta \int_{-k}^k \phi(w-ae) [L^{(k+\varepsilon)}_\beta (w)-L^{(k)}_\beta (w)] dw \\
& \hskip 2em + 2 \beta \phi(k-ae)L^{(k+\varepsilon)}_\beta (k+\varepsilon)\varepsilon + O(\varepsilon^2)
\end{align*}

Let $R^{(k)}_\beta(e,w;a)$ be the resolvent of $\phi$, as given in~(16). Then,
\begin{align*}
L^{(k+\varepsilon)}_\beta (e) - L^{(k)}_\beta (e) &=2 \beta \int_{-k}^k R^{(k)}_\beta(e,w;a)\phi(k-ae)L^{(k+\varepsilon)}_\beta (w)\varepsilon dw \\
&\hskip 4em +  O(\varepsilon^2)
\end{align*}
This implies that
\begin{align*}
\Bigm|\frac{L^{(k+\varepsilon)}_\beta (e)-L^{(k)}_\beta (e)}{\varepsilon}\Bigm|
&\le 2 \|\phi\|_\infty \|L^{(k)}_\beta\|_\infty \Bigm| \int_{-k}^k \beta R^{(k)}_\beta(e,w;a) dw\Bigm|\\
 &\hskip 3em + O(\varepsilon).
\end{align*}
Since $\beta \mathcal{B}^{(k)}$ is a contraction, the value of the integral in the first term on the right hand side of the above inequality is less than 1 and the result follows from the definition of differtiability.

\end{document}